\documentclass{amsart}

\usepackage{pb-diagram}
\usepackage{amssymb}

\theoremstyle{plain}
\newtheorem{theorem}{Theorem}[section]
\newtheorem{corollary}[theorem]{Corollary} 
\newtheorem{lemma}[theorem]{Lemma}
\newtheorem{proposition}[theorem]{Proposition}

\newtheorem*{theorema}{Theorem A}
\newtheorem*{theoremb}{Theorem B}

\numberwithin{equation}{section}

\theoremstyle{definition}
\newtheorem{definition}[theorem]{Definition}
\newtheorem{example}[theorem]{Example}

\theoremstyle{remark}
\newtheorem{remark}[theorem]{Remark}

\newcommand{\la}{\langle}
\newcommand{\ra}{\rangle}

\newcommand{\ow}{\omega}

\newcommand{\R}{\mathbb{R}}
\newcommand{\C}{\mathbb{C}}
\newcommand{\Z}{\mathbb{Z}}
\newcommand{\N}{\mathbb{N}}
\newcommand{\orbit}{\mathcal{O}}

\newcommand{\fg}{\mathfrak{g}}
\newcommand{\fh}{\mathfrak{h}}

\newcommand{\Su}{\operatorname{SU}}
\newcommand{\So}{\operatorname{SO}}

\newcommand{\pz}{\partial z}

\newcommand     {\comment}[1]   {}
\newcommand{\mute}[2] {}
\newcommand     {\printname}[1] {}

\newcommand{\labell}[1] {\label{#1}\printname{#1}}

\newcommand{\Yone}{\ensuremath{\So(3)\times_{S^1}\R^2\times\C}}
\newcommand{\Ytwo}{\ensuremath{\So(3)\times_{N(S^1)}\R^2\times\C}}
\newcommand{\Ythree}{\ensuremath{\So(3)\times_{N(S^1)}(\R^2\times\C)}}

\begin{document}

\title{Complexity one Hamiltonian $\Su(2)$ and $\So(3)$ actions}
\author{River Chiang}
\begin{abstract}
We consider compact connected six dimensional symplectic manifolds
with Hamiltonian $\Su(2)$ or $\So(3)$ actions with cyclic principal 
stabilizers. 
We classify such manifolds up to equivariant symplectomorphisms. 
\end{abstract}
\maketitle

\tableofcontents

\comment{$\overline{\Phi}$ is induced by the norm square of the moment map}

\comment{Need to check D-H}

\comment{$G$-equivariant, $G$-invariant, 4-manifold, 6 dimensional, 
$G$-mfld, $G$-action, $SO(3)$-mfold, $SO(3)$ action}

\section{Introduction}\labell{intro}

Let $(M, \ow)$ be a symplectic manifold and $G$ be a compact
connected Lie group that acts effectively on $M$ by 
symplectic transformations. 
A \textbf{moment map}
$\Phi\colon M \to \fg^*$ is a $G$-equivariant map such that for every 
$\xi$ in the Lie algebra $\fg$ of $G$, 
\[
\iota(\xi_M)\ow = - d\la \Phi, \xi\ra
\]
where $\xi_M \colon M \to TM$ denotes the induced vector field of $\xi$
on $M$. If there is a moment map, we say that the action is
\textbf{Hamiltonian}. 
The triple $(M, \ow, \Phi)$ is called a \textbf{Hamiltonian
$\mathbf{G}$-manifold} or \textbf{Hamiltonian $\mathbf{G}$-action}.
An isomorphism between two such manifolds is an equivariant
symplectomorphism that respects the moment maps. We usually assume that 
$M$ is connected, that $G$ acts effectively on $M$, and that the moment
map $\Phi$ is proper.

For a point $\alpha$ in the dual of the Lie algebra of $G$, the 
\textbf{symplectic quotient} or \textbf{reduced space} at $\alpha$ is
the topological space $M_\alpha =  
\Phi^{-1}(G \cdot \alpha) / G = \Phi^{-1}(\alpha) / G_\alpha$,
where $G \cdot \alpha$ is the coadjoint orbit through $\alpha$ and 
$G_\alpha$ denotes the stabilizer of $\alpha$. 
If $\alpha$ is a regular value of the moment 
map $\Phi$, the 
reduced space $M_\alpha$ 
is a symplectic orbifold. 
In general, the reduced space is a symplectic stratified space
\cite{ls, bl}.
The \textbf{complexity} of $(M, \ow, \Phi)$ is half the 
dimension of the reduced space $M_\alpha$ at a generic
value $\alpha$ in the moment image $\Phi(M)$.  

Suppose a torus $T$ acts on a symplectic manifold $M$ in a
Hamiltonian fashion. Its complexity is $\frac{1}{2}\text{dim}M -
\text{dim}T$. In particular, the complexity is zero exactly if the
dimension of the torus is half the dimension of the manifold. The
space is then called a toric manifold. These spaces are classified
by their moment images \cite{dE1}. Karshon and Tolman have studied
complexity one torus actions in arbitrary dimensions \cite{kt1, kt2}.
Special cases are compact symplectic four-manifolds with circle
actions which are classified by Karshon \cite{k1}, also see Audin,
and Ahara-Hattori \cite{aU1, aU2, ah}.
 
As for $\So(3)$ actions, there are two distinct cases: the principal
stabilizer is $S^1$ or the principal stabilizer is $\Z_n = \Z/n\Z$. 
The former case is characterized by Iglesias \cite{iG}. 
In any dimension, the manifold $M$ is isomorphic to the product of 
$S^2$ by the symplectic orbit manifold $M/\So(3)$. 
The latter case in dimension four has complexity zero and is  
classified by Iglesias \cite{iG}. The only compact 
Hamiltonian $\So(3)$-manifolds are $\C P^2$ and $S^2 \times
S^2$. The first is equipped with the natural
action induced by $\Su(3)$ and the second can be equipped with
different $\So(3)$ actions indexed by $\N$.
Complexity zero Hamiltonian actions
of more general nonabelian groups have been studied by Delzant, 
Woodward, and Knop \cite{dE2, w, kN}. 

In the algebraic and smooth categories, Lie group actions of complexity
zero or one have been studied in \cite{tI1, tI2, fI, ow}.

In this paper we study complexity one $\Su(2)$ and $\So(3)$
actions. After Iglesias's work, it remains to classify
compact connected six dimensional symplectic
manifolds equipped with Hamiltonian $\Su(2)$ or $\So(3)$ 
actions with cyclic principal stabilizers  
$\Z_n = \Z/n\Z$.  
When the moment image does not contain zero, the manifold is of the
form $G \times_{S^1} X$ where $G$ is $\Su(2)$ or $\So(3)$ and $X$ is
a symplectic four-manifold with a (possibly noneffective) Hamiltonian 
circle action. This can be viewed as an 
immediate corollary from the classification of circle actions on
four-manifolds. Therefore, we emphasize the case when zero is in
the moment image.  

Generalizing techniques established by Karshon and
Tolman \cite{kt1, kt2}, we proceed by first studying
the basic building blocks: the preimages
under the moment map of sufficiently small open subsets in $\fg^*$.
In this context, the classification applies not
only to a compact manifold but also to a noncompact manifold with a
\textbf{proper} moment map. 
In Sections \ref{S:part1begin}--\ref{S:part1end}, 
we provide a complete set of invariants for the
preimage of a neighborhood of $0 \in \fg^*$. 
In Sections \ref{S:part2begin}--\ref{S:part2end}, 
we discuss local invariants for the preimage of any neighborhood
in $\fg^*$ away from zero. We then show that if two spaces are
locally isomorphic, they are globally isomorphic.

We now describe the invariants.

The \textbf{Duistermaat-Heckman function} is a real function defined
on the dual of the Lie algebra that takes the value of the symplectic
volume of the reduced space. That is, $f\colon\fg^* \to \R$ such that 
\[
f(\alpha)=\text{Vol}(M_\alpha).
\] 

\mute{old definition}{
%The \textbf{Liouville measure} on a $2n$ dimensional symplectic
%manifold $(M, \ow)$ is given by integration of the volume form
%$\ow^n/n!$ with respect to the symplectic orientation. In the
%presence of a Hamiltonian $G$ action, the
%\textbf{Duistermaat-Heckman measure} is the push-forward of Liouville
%measure by the moment map. It is equal to Lebesgue measure on $\fg^*$
%times the \textbf{Duistermaat-Heckman function}.
}

For $x \in M$, the set $G_x = \{\, g \in G \, | \, gx=x\,\}$ of
elements of $G$ leaving $x$ fixed is called the \textbf{stabilizer} 
of $x$. $G_x$ is a closed subgroup of $G$. It acts linearly on
the tangent space $T_xM$, and thus $T_xM$ is a representation space
of $G_x$.  This representation is called the
\textbf{isotropy representation} at $x$. The stabilizers of points in
the same orbit are conjugate, and 
their isotropy representations are linearly symplectically
isomorphic. We call this conjugacy class the \textbf{stabilizer},
and the isomorphism class the \textbf{isotropy
representation} of the orbit.  

For convenience, we call the level set of the moment map 
$\Phi^{-1}(\alpha)$ the \textbf{moment fiber} at $\alpha$ or simply the 
$\mathbf{\alpha}$ \textbf{fiber}.  
An orbit is \textbf{exceptional} if it has a
strictly larger stabilizer than any nearby orbit in the same moment
fiber.  
In particular, if $\Phi^{-1}(\alpha)$ contains only one
$G_\alpha$ orbit, that orbit is exceptional. Since the moment map is proper, 
each moment fiber is compact, and it has finitely many exceptional
orbits.  The \textbf{isotropy data} at $\alpha$ consist of 
the unordered list
of isotropy representations of the exceptional orbits in
$\Phi^{-1}(\alpha)$.

\comment{level set = $\Phi^{-1}(G\cdot\alpha)$ or
$\Phi^{-1}(\alpha)$? fiber vs. level set?}

\comment{I kind of feel bad that the intro jumps to $G_\alpha$ orbit
instead of $G$-orbit. maybe 
just say that there is a one-to-one correspondence between $G$-orbit and 
$G_\alpha$-orbit? or is this too well-known to mention?}

If $M$ is a $G$-manifold with connected orbit space $M/G$, there
exists an open dense subset of $M$ in which all stabilizers are
conjugate. This conjugacy class is called the 
\textbf{principal stabilizer} of $M$. A
similar notion exists for the zero fiber of the moment map.

If $\Phi^{-1}(\alpha)$ consists of a single $G_\alpha$-orbit, it is
called a \textbf{short} fiber; otherwise, it is \textbf{tall}.
In the complexity one case, if $\Phi^{-1}(\alpha)$ is tall,
we show that its reduced space 
$M_\alpha = \Phi^{-1}(\alpha)/G_\alpha$ is topologically a closed 
connected oriented surface. 
We call its genus the \textbf{genus at} $\mathbf{\alpha}$. 
Defining the genus of a point to be zero, we show that the genus is
independent of $\alpha$ for any $\alpha \in \Phi(M)$ and is called 
the \textbf{genus} of the Hamiltonian $G$-manifold $(M, \ow, \Phi)$. 

For any real $n$ dimensional vector bundle $\pi\colon W \to X$, 
there exists an associated orientation bundle
$p\colon \Tilde{X} \to X$, whose fiber over a point $x$ is the 
two ways to orient $\pi^{-1}(x)$. This is a two-sheeted covering, and 
\v{C}ech cocycles provide a convenient way to construct it. 
Choose an open cover $\mathcal{U} = \{U_\alpha\}$ of $X$ with
trivialization maps 
$\varphi_\alpha\colon  U_\alpha \times \R^n \to \pi^{-1}(U_\alpha)$. 
The Jacobian determinants of the change of fiber coordinates from 
$\R^n$ to $\R^n$ have a locally constant sign, 
which gives a locally constant function from 
$U_\alpha \cap U_\beta$ to $\Z_2$. The chain rule for Jacobians
implies that this is a cocycle. 
%Define $p \colon  \Tilde{X} \to X$ to be
%the resulting $\Z_2$ bundle. This is the orientation bundle
%associated with the vector bundle $W$ and 
It determines
an element $w_1(W) \in H^1(X;\Z_2)$, called the \textbf{first
Stiefel-Whitney class} of $W$. In a similar fashion, 
we can construct the associated orientation bundle and 
define the first Stiefel-Whitney class of any fiber bundle $W$ 
over a manifold $X$ provided that the fiber of $W$ is 
connected and orientable. 

In particular, when the zero fiber of a Hamiltonian $\So(3)$-manifold
is tall, and when the principal stabilizer of the zero fiber 
is $S^1$, the zero fiber $\Phi^{-1}(0)$ is
a sphere bundle over the reduced space 
off the exceptional orbits. Let $E = \{\,E_j\,\}$ denote the set of
exceptional orbits in the zero fiber.
Let $M^{\text{reg}}_0$ denote the smooth part of the symplectic
quotient at $0$, i.e., $M_0^{\text{reg}}=(\Phi^{-1}(0) \smallsetminus
E)/G$. Then $M_0^{\text{reg}}$ is diffeomorphic to 
$\Sigma \smallsetminus \{\text{\emph{finitely many points}}\}$,
where $\Sigma = \Phi^{-1}(0)/G$ is a closed connected oriented
surface.  
Through the orientations on the fiber spheres,
$\Phi^{-1}(0)$ 
induces an associated orientation bundle on $M_0^{\text{reg}}$,
and the first Stiefel-Whitney class in $H^1(M_0^{\text{reg}};\Z_2)$.

Sections \ref{S:part1begin}--\ref{S:part1end} 
are devoted to prove the local uniqueness over $0$:

\begin{theorema}[Local Uniqueness over $0$]\labell{T:locun}
Let $G$ be $\Su(2)$ or $\So(3)$. 
Let $(M,\ow, \Phi)$ and $(M', \ow', \Phi')$ be 
compact connected six dimensional Hamiltonian $G$-manifolds such that 
$0 \in \Phi(M) = \Phi'(M')$. 
There exists an invariant neighborhood $V$ of $0$ in $\fg^*$ over which the 
Hamiltonian $G$-manifolds are isomorphic if and only if 
\begin{itemize}
\item their Duistermaat-Heckman functions coincide on $V$, 
%(in $\So(3)$ case, they only have to agree at some $\alpha \in V$), 
\item their isotropy data and genus at $0$ are the same,
\item their principal stabilizers of the zero fibers are the same,
\item if the zero fibers are tall with principal stabilizer $S^1$, 
the first Stiefel-Whitney classes of $\Phi^{-1}(0)$ and
${\Phi'}^{-1}(0)$ in $H^1(M_0^{\text{reg}}; \Z_2)$ and
$H^1({M'}_0^{\text{reg}}; \Z_2)$
are equal (under an identification of $M_0$ and $M_0'$ that respects
the isotropy data). 
\end{itemize}
\end{theorema}

\begin{remark}
For $G=\Su(2)$, a tall zero fiber has no exceptional orbits. In this
case, Theorem A follows from the equivariant symplectic
embedding theorem of \cite{w1}. 
\end{remark}

In Sections \ref{S:part2begin}-\ref{S:part2end}, 
we adapt the idea of symplectic cross-sections introduced by Guillemin
and Sternberg \cite{gs2}. It allows us to apply our techniques 
and determine when two Hamiltonian $G$-manifolds are locally isomorphic.
We then define compatible invariants to construct a global
isomorphism from the local isomorphisms. 

Let $E$ denote the set of exceptional orbits in $M$. We consider the
projections $M\to M/G$ and $\fg^* \to \fg^*/G$, and the 
map $\overline{\Phi}$ induced by the moment map $\Phi$.
The 
\textbf{isotropy skeleton} is the space $E/G$ where each point is 
labeled by its isotropy representation, together with the map 
$\overline{\Phi}\colon  E/G \to \fg^*/G$. Two isotropy skeletons are 
considered the same if there exists a homeomorphism $f\colon  E/G \to E'/G$
that sends each point to a point with the same isotropy representation
and such that $\overline{\Phi} = \overline{\Phi'}\circ f$. 

\comment{Is it OK that I use $E$ for the set of the exceptional
orbits in $\Phi^{-1}(0)$ before but use it now for those in $M$?}

We have the following global uniqueness theorem:

\begin{theoremb}\labell{T:globun}
Let $G$ be $\Su(2)$ or $\So(3)$. 
Let $(M, \ow, \Phi)$ and $(M', \ow',\Phi')$ be compact connected
six dimensional Hamiltonian $G$-manifolds such that 
$\Phi(M)=\Phi'(M)$. Then $M, M'$ are isomorphic if and only if
they have the same Duistermaat-Heckman function,
the same genus, the same isotropy skeleton, the same principal 
stabilizers of the manifolds and of the zero fibers,  
and the same first Stiefel-Whitney class of the zero fibers when
applicable.  
\end{theoremb}

\section{The zero fiber}\labell{S:part1begin}

We begin by stating some properties of the zero fiber of the
moment map. 
If a compact Lie group acts on a symplectic manifold and if the
action is Hamiltonian, the Local Normal Form Theorem of Marle, Guillemin
and Sternberg \cite{mA, gs1} provides a nice description of the 
neighborhood of any orbit in the zero fiber. 

\begin{theorem}[Local Normal Form]\labell{T:normalform}Let a compact Lie group
$G$ act on a symplectic manifold $(M, \ow)$ with a moment map $\Phi\colon 
M \to \fg^*$. Let $x$ be a point in the zero fiber of the moment map,
$G \cdot x$ be its orbit in $M$, $H$ be the stabilizer of $x$, 
and $V$ be the symplectic slice $(T_x(G\cdot x))^\ow / T_x(G \cdot
x)$ at $x \in M$.
Given a choice of an $H$-equivariant
splitting $\fg = \fh \oplus \mathfrak{m}$, there exists a
$G$-invariant symplectic form on the \textbf{local model}
$Y = G \times_{H} (\fh^0 \times V)$ such that 
\begin{enumerate}
\item a neighborhood of $G \cdot x$ is equivariantly
symplectomorphic to a neighborhood of the zero section in $Y$, and
\item the action of $G$ on $Y$ is Hamiltonian and the moment map is
given by 
\[
\Phi_Y([g, \mu, v]) = Ad^\dag (g)(\mu + \pi^*\Phi_V(v)),
\]
\end{enumerate}
where $\fh^0$ is the annihilator of $\fh$,
$Ad^\dag$ is the coadjoint action, $\pi^*\colon  \fh^* \to \fg^*$
is induced by the projection $\pi\colon  \fg \to \fh$,  
and $\Phi_V\colon  V \to \fh^*$ is the moment map
for the slice representation. 
\end{theorem}

A special case of Theorem A follows immediately
from the Local Normal Form Theorem: 

\begin{proposition}\labell{P:shortzerofiber}
Let $G=\Su(2)$ or $\So(3)$. 
Let $(M, \ow, \Phi)$ and $(M', \ow', \Phi')$ be Hamiltonian 
$G$-manifolds such that $0 \in \Phi(M) = \Phi'(M')$. 
Assume $\Phi^{-1}(0)$ and ${\Phi'}^{-1}(0)$ 
consist of one single orbit each and that these 
orbits have the same isotropy 
representation. Then there exists a neighborhood $V$ of $0$ in
$\fg^*$ over which $M$ and $M'$ are isomorphic. 
\end{proposition}

\begin{proof}
Since $\Phi^{-1}(0)$ and ${\Phi'}^{-1}(0)$ 
consist of one single orbit each with the same 
isotropy representation, we can find $x \in M$ and $x' \in M'$ with
the same stabilizer such that $G \cdot x = \Phi^{-1}(0)$ and $G \cdot
x' = {\Phi'}^{-1}(0)$.   
It follows from the Local Normal Form Theorem that there exist
neighborhoods $U$ of $G \cdot x$ and $U'$ of $G \cdot x'$ and an
equivariant symplectomorphism $\varphi\colon U \to U'$ such that
$\varphi(x)=x'$ and $\Phi' \circ \varphi = \Phi$. 

Since the moment maps $\Phi$ and $\Phi'$ are proper, 
there exist neighborhoods $W$ and $W'$ of $0$ in $\fg^*$ 
such that $\Phi^{-1}(W) \subset U$ and 
${\Phi'}^{-1}(W') \subset U'$. 
We can then take $V = W \cap W'$. 
\end{proof}

If we ought to understand all the possible local models 
$G \times_{H} (\fh^0 \times V)$, we first have to 
understand all the isotropy representations.  
The isotropy representation is a
direct sum of the coadjoint action of $H \subset G$ on $\fh^0 \subset
\fg^*$ and the slice representation of $H$ on $V$. Therefore, we
need to know all the possible stabilizers, 
slice representations, and coadjoint actions for $\Su(2)$ 
and $\So(3)$.  

Up to conjugacy, the finite subgroups of $\So(3)$ include the trivial
group $\{1\}$, the cyclic groups $Z_k, \, k=2, 3, \dots$, the dihedral
groups $D_{2k}, \, k=2, 3, \dots$, the tetrahedral, octahedral 
and icosahedral groups. Any of
these finite subgroups will be denoted by $\Gamma$ if it doesn't need
to be specified.

\comment{cite?}

Up to conjugacy, $\So(3)$ has two infinite one dimensional closed 
subgroups: the maximal abelian subgroup $S^1$, and its normalizer, 
which is isomorphic to $\operatorname{O}(2)$ and will be denoted as
$N_{\So(3)}(S^1)$, or simply $N(S^1)$ if there is no possible confusion. 

Stabilizers are closed subgroups. Therefore, the possible stabilizers of an
$\So(3)$ action, up to conjugacy, 
are the subgroups listed above and $\So(3)$ itself. 

Similarly, up to conjugacy, the closed subgroups of $\Su(2)$ 
are: a collection of finite subgroups,  
again denoted by $\Gamma$, the maximal torus $S^1$, 
the normalizer of the maximal torus $N_{\Su(2)}(S^1)$, 
and $\Su(2)$ itself.  
Note that $N_{\Su(2)}(S^1)$ is no longer isomorphic 
to $\operatorname{O}(2)$.

The slice representations we need to consider are linear symplectic
representations of $H$ on $\C$ and those of 
$G$ on $\C^3$, where $G$ denotes either $\Su(2)$ or $\So(3)$ and $H$
denotes $S^1$ or $N_G(S^1)$. The linear symplectic representations
on a complex vector space $\C^n$ 
are equivalent to the unitary representations. Therefore, we know that
the representations of $S^1$ on $\C$ are
characterized by the weights $n \in \Z$, and that there is only one
effective representation each for $\Su(2)$ and
$\So(3)$ on $\C^3$. 

A slice representation 
$\rho\colon  N_G(S^1) \to \operatorname{U}(\C) = S^1$
is an analytic homomorphism. 
Because $S^1$ is abelian, the 
commutator group of $N_G(S^1)$ is in the kernel of $\rho$. 
So the kernel of $\rho$ is either $N_G(S^1)$, or $S^1$. 
The former implies that the slice representation is trivial; the
latter implies that the slice representation reduces 
to a $\Z_2$ action such that $h \cdot z = z$ for $h \in S^1$ and $h
\cdot z = -z$ otherwise.

We fix an inner product on the Lie algebra $\fg$ of $G = \Su(2)$
or $\So(3)$. This determines a projection $\fg \to \fh$ and 
the induced inclusion $\fh^* \to \fg^*$ for any $\fh \subseteq \fg$. We
also identify $\fh^*$ with its image in $\fg^*$. 
Since maximal tori in the same group are conjugate to each
other, we use 
$\left(\begin{smallmatrix}
e^{i\theta} & 0 \\
0 & e^{-i\theta}
\end{smallmatrix}\right)$ to represent $S^1$ in $\Su(2)$ and 
$\left(\begin{smallmatrix}
1 & 0 & 0 \\
0 & \cos\theta & -\sin\theta \\
0 & \sin\theta & \cos\theta 
\end{smallmatrix}\right)$ in $\So(3)$ for the standard local models.
We fix the identification between $\fg^*$ and 
$\R^3$ throughout this paper so that the annihilator of the above chosen 
$S^1$ is identified with 
$\left\{\left(\begin{smallmatrix}0\\x\\y\end{smallmatrix}\right)\right\}
\simeq \R^2 \subset \R^3$. 
With these identifications, the coadjoint actions of $\So(3)$ and
$\Su(2)$ are rotations at speed 1 and 2, such that $\fg^*$ is fixed 
by the centers $\operatorname{I}$ and $\Z_2 = \pm \operatorname{I}$.

\comment{check consistency throughout the paper about $\R^2=(0, x, y)$}

\begin{corollary} Let $(M, \ow, \Phi)$ be a six dimensional
Hamiltonian $\Su(2)$-manifold with a moment map 
$\Phi\colon M \to \mathfrak{su}(2)^* \simeq \R^3$. 
Assume the action is effective. Then for any
$x \in \Phi^{-1}(0)$, the local model for the orbit $G\cdot x$ is one
of the following:
	\begin{enumerate}
	\item $Y = \Su(2) \times_{\Gamma} \R^3 = \{\,[g,
	\mu]\,|\,g\in \Su(2), \mu \in \R^3, [g,\mu]
	=[ga^{-1}, Ad^\dag(a)\mu],
	\forall a \in\Gamma\,\}$ with 
	$\Phi_Y([g, \mu])= Ad^\dag (g)\mu$, where $\Gamma$ does not contain 
	$\Z_2 = \pm \operatorname{I}$. 
	\item $Y = \Su(2) \times_{S^1} (\R^2 \times \C)$, where $S^1$ acts
	on $\C$ with an odd numbered weight $n$, and $\Phi_Y([g, \mu, z]) 
	= Ad^\dag (g)(\mu + \frac{n}{2}|z|^2)$.
	\item $Y = \C^2 \times \C = \{\,(u, v, w)\, | \, u, v, w \in
	\C\,\}$,  where $\Su(2)$ acts on $\C^2$ as  
	the standard unitary transformations
	of $\C^2$ to itself. The moment map on $Y$ is 
	$\Phi(u,v,w) = \left(\frac{|u|^2-|v|^2}{2},
	\text{Re}\,(\bar{u}v), \text{Im}\,(\bar{u}v)\right)$.
	\end{enumerate}
\end{corollary}

\comment{better way to describe local models? 
still need to mention effective in the statements?}

\comment{$X_1=\left(\begin{smallmatrix}i&0\\0&-i\end{smallmatrix}\right)$, 
$X_2=\left(\begin{smallmatrix}0&i\\i&0\end{smallmatrix}\right)$, 
$X_3=\left(\begin{smallmatrix}0&1\\-1&0\end{smallmatrix}\right)$, 
doesn't really matter that much}

\begin{corollary} Let $(M, \ow, \Phi)$ be a six dimensional Hamiltonian
$\So(3)$-manifold with a moment map $\Phi\colon 
M \to \mathfrak{so}(3)^* \simeq \R^3$. Assume the action is effective.
Then for any 
$x \in \Phi^{-1}(0)$, the local model for the orbit $G\cdot x$ is one
of the following:
	\begin{enumerate}
	\item $Y = \So(3) \times_{\Gamma} \R^3$, and 
	$\Phi_Y([g, \mu])= Ad^\dag (g)\mu$. 
	\item $Y = \So(3) \times_{S^1} (\R^2 \times \C)$, where $S^1$ acts
	on $\C$ with weight $n$, and $\Phi_Y([g, \mu, z]) 
	= Ad^\dag (g)(\mu + \frac{n}{2}|z|^2)$.
	\item $Y = \So(3) \times_{S^1} \R^2 \times \C$, where $S^1$
	acts trivially on $\C$, and $\Phi_Y([g, \mu], z) 
	= Ad^\dag (g)\mu$. 
	\item $Y = \So(3) \times_{N(S^1)} (\R^2 \times \C)$, where
	$N(S^1)$ acts on $\C$ as $N(S^1)/S^1 \simeq \Z_2$, and 
	$\Phi_Y([g, \mu, z])=Ad^\dag(g)\mu$.
	\item $Y = \So(3) \times_{N(S^1)} \R^2 \times \C$, where
	$N(S^1)$ acts trivially on $\C$, and 
	$\Phi_Y([g, \mu], z)=Ad^\dag(g)\mu$.
	\item $Y = \C^3 = \{\, q + \sqrt{-1}\,p \, |\, q, p \in \R^3\} 
	\simeq T^*\R^3$, where $\So(3)$ acts on $\C^3 = T^*\R$ through
	its standard action on $\R^3$, and the moment map is the vector 
	cross product,
	i.e. $\Phi_Y(q, p) = q \times p$.
	\end{enumerate}
\end{corollary}

\begin{remark}
For convenience, 
from now on, we will refer to the local models using the expressions
appeared in the above corollaries. In particular, the isotropy
representations are implied when we use different expressions. 
For instance, the local model $\C^2 \times \C$ 
has an $\Su(2)$ action on the first $\C^2$ while 
the local model $\C^3$ has a
canonical $\So(3)$ action. 
In addition, $G \times_H (\R^2 \times \C)$
always has a nontrivial action of $H$ on $\C$ for $G=\Su(2)$ or
$\So(3)$ and $H = S^1$ or $N(S^1)$.
\end{remark}

We can read out information from the local models. 
For example, since any moment fiber is connected, and 
since any orbit is closed, 
if there is only one orbit $\orbit \in \Phi^{-1}(0)$
sitting inside a local model, there is
only one orbit in $\Phi^{-1}(0)$ and therefore 
$\Phi^{-1}(0)$ is a short fiber.

\begin{corollary} Let $G$ be $\Su(2)$ or $\So(3)$.
Let $(M, \ow, \Phi)$ be a six dimensional 
Hamiltonian $G$-manifold. Assume the zero fiber is short. The local
model for the single orbit in the zero fiber is either
$G \times_{\Gamma} \R^3$, or $G \times_{S^1} (\R^2 \times \C)$.
\end{corollary}

\mute{muting the definition of exceptional orbit}{
\begin{definition}
We call an orbit $\orbit$ an \textbf{exceptional orbit} if it has 
strictly larger stabilizer than any nearby orbit in the same moment
fiber.  In particular, if $\Phi^{-1}(\alpha)$ contains only one
$G_{\alpha}$ orbit,  that orbit is exceptional. 
\end{definition}}

\begin{corollary}
Let $G$ be $\Su(2)$ or $\So(3)$.
Let $(M, \ow, \Phi)$ be a six dimensional 
Hamiltonian $G$-manifold. Assume an orbit $\orbit$ in $\Phi^{-1}(0)$ is 
nonexceptional. Then the local model for $\orbit$ is either 
$\So(3) \times_H \R^2 \times \C$ with $H = S^1$ or 
$N_{\So(3)}(S^1)$, or $\C^2 \times \C$. 
\end{corollary}

The following theorem from Lerman and Sjamaar \cite{ls} 
states that the reduced space at $0$ is stratified:

\begin{theorem}\labell{T:strata}
Let $G$ act on a symplectic manifold $M$ with a proper moment map.
The reduced space $M_0 = \Phi^{-1}(0)/G$ can be decomposed into a 
disjoint union of symplectic manifolds with respect to the stabilizers,
i.e. \
\[
M_0 = \bigsqcup_{H < G} (M_0)_H
\]
where $(M_0)_H$ denotes all the orbits in the zero fiber 
whose stabilizer is conjugate to $H$. 

Moreover, there exists a unique piece  
$(M_0)_H$ which is open,
connected, and dense in the reduced space.
\end{theorem}

\begin{definition}
The stabilizer $H$ of the unique piece $(M_0)_H$ 
stated in Theorem \ref{T:strata} is 
called the \textbf{principal stabilizer} of the zero fiber. 
\end{definition}

\begin{corollary}
Let $G$ be $\Su(2)$ or $\So(3)$.
Let $(M, \ow, \Phi)$ be a six dimensional 
Hamiltonian $G$-manifold. Assume 
that the zero fiber $\Phi^{-1}(0)$ is tall.
Then the principal stabilizer of the zero fiber is either
\begin{itemize}
\item $\Su(2)$, when $G = \Su(2)$, or
\item $S^1$ or $N_G(S^1)$, when $G = \So(3)$.
\end{itemize}
\end{corollary}

\begin{corollary}
Let $G$ be $\Su(2)$ or $\So(3)$.
Let $(M, \ow, \Phi)$ be a six dimensional 
Hamiltonian $G$-manifold. Assume  
the zero fiber $\Phi^{-1}(0)$ is tall. Then
\begin{itemize}
\item there is no exceptional orbit in the zero fiber
when $G = \Su(2)$;
\item 
when $ G = \So(3)$, if there exists an exceptional orbit in the zero 
fiber, the local model for the exceptional orbit is either
$\So(3) \times_{N_G(S^1)} \R^2 \times \C$ or
$\C^3$.  In this case, the principal stabilizer of the zero fiber
is $S^1$.  
\end{itemize}
\end{corollary}

The list of local models and their properties 
is important and useful later on. We 
use Table \ref{Ta:model} for an easy reference. 

\begin{table}[h]
\centering
\begin{tabular}{|l|l|l|l|l|}\hline
\multicolumn{5}{|c|}{$\Su(2)$}\\ \hline \hline
\emph{local model } & 
\emph{moment map} & 
\emph{zero} &
\emph{central orbit } &
\emph{principal} \\ 
$G \times_{H} (\fh^0 \times V)$ & $\Phi$ & \emph{fiber} 
& $\{[g, 0, 0]\}$ & \emph{stabilizer}  \\ 
& & & & \emph{of} $\Phi^{-1}(0)$\\ \hline
$\Su(2)\times_\Gamma \R^3$ & 
$Ad^\dag(g)\mu$ &
short &
exceptional &
$\Gamma$ \\ \hline
$\Su(2)\times_{S^1} (\R^2 \times \C)$ & 
$Ad^\dag(g)(\mu + \frac{n}{2}|z|^2)$ &
short &
exceptional &
$S^1$ \\ \hline
$\C^2 \times \C$ & 
$(\frac{|u|^2-|v|^2}{2}, \bar{u}v)$ &
tall &
nonexceptional &
$\Su(2)$ \\ \hline \hline
\multicolumn{5}{|c|}{$\So(3)$}\\ \hline \hline
\emph{local model } & 
\emph{moment map} & 
\emph{zero} &
\emph{central orbit } &
\emph{principal} \\ 
$G \times_{H} (\fh^0 \times V)$ & $\Phi$ & \emph{fiber} 
& $\{[g, 0, 0]\}$ & \emph{stabilizer}  \\ 
& & & & \emph{of} $\Phi^{-1}(0)$\\ \hline
$\So(3)\times_\Gamma \R^3$ & 
$Ad^\dag(g)\mu$ &
short &
exceptional &
$\Gamma$ \\ \hline
$\So(3)\times_{S^1} (\R^2 \times \C)$ & 
$Ad^\dag(g)(\mu + \frac{n}{2}|z|^2)$ &
short &
exceptional &
$S^1$ \\ \hline
$\So(3)\times_{N(S^1)} \R^2 \times \C$ & 
$Ad^\dag(g)\mu$ &
tall &
nonexceptional &
$N(S^1)$ \\ \hline
$\So(3)\times_{S^1} \R^2 \times \C$ & 
$Ad^\dag(g)\mu$ &
tall &
nonexceptional &
$S^1$ \\ \hline
$\So(3)\times_{N(S^1)} (\R^2 \times \C)$ & 
$Ad^\dag(g)\mu$ &
tall &
exceptional &
$S^1$ \\ \hline
$\C^3$ & 
$q \times p$ &
tall &
exceptional &
$S^1$ \\ \hline 
\end{tabular}\vspace{.2cm}
\caption{A list of local models.}\labell{Ta:model}
\end{table}

\section{Eliminating the symplectic form}\labell{S:eliminatew}

An equivariant symplectomorphism between two Hamiltonian 
$G$-manifolds preserves the orientations. Using Moser's method
\cite{w1}, we show 
that under appropriate circumstances, we can recover the symplectic form
from an orientation preserving equivariant
diffeomorphism that respects the moment maps. This enables us to use
the techniques in differential topology in later sections. 

\begin{definition} \labell{D:phiGdiffeo}
Let a compact Lie group $G$ act on oriented manifolds
$M$ and $M'$ with 
$G$-equivariant maps $\Phi\colon M \to \fg^*$ and $\Phi'\colon M' \to \fg^*$.
A \textbf{$\mathbf{\Phi}$-$\mathbf{G}$-diffeomorphism} 
$\Psi\colon M\to M'$ is an orientation  
preserving equivariant diffeomorphism such that 
$\Psi^*(\Phi') = \Phi$.
\end{definition}

Here we assume a technical condition:
\begin{equation}\labell{E:tech}
\begin{split}
&\text{The restriction map } H^i(N/G) \to 
H^i(\Phi^{-1}(G \cdot \alpha)/G), \quad i=1, 2 \\
&\text{is one-to-one with $\Z$ or $\Z_2$ coefficient}\\
&\text{for some neighborhood $N$ of the zero fiber 
and } \forall\, \alpha \in \Phi(N) \subset \fg^*.    
\end{split}
\end{equation}
We will prove later in this paper that such a neighborhood $N$ always
exists and this condition is satisfied. 

\begin{remark}
Through out this paper, when we mention a neighborhood $N$ of the
zero fiber, we always assume that $N$ 
is small enough so that it can be covered by invariant
open subsets of the local models in Table \ref{Ta:model}.  
\end{remark}

\begin{proposition}\labell{P:symplecto} Let $G$ be $\Su(2)$ or
$\So(3)$. 
Let $(M, \ow, \Phi)$ and $(M', \ow', \Phi')$ be 
neighborhoods of the zero fibers in six dimensional Hamiltonian
$G$-manifolds such that $0 \in \Phi(M) = \Phi(M')$.  
Assume that $M$ and $M'$ satisfy Condition \ref{E:tech} and 
that their Duistermaat-Heckman functions coincide.
Then there exist neighborhoods $N \subset M$ and $N'
\subset M'$ of the zero fibers   
such that there exists an equivariant symplectomorphism from
$N$ to $N'$ 
if and only if there exists a $\Phi$-$G$-diffeomorphism from 
$N$ to $N'$.
\end{proposition}

We prove this with two lemmas. First we
recall that a differential form $\beta$ on $M$ is \textbf{basic} if
it is $G$-invariant and $\iota(\xi_M)\beta =0$ 
for every $\xi \in \fg$ and $\xi_M$ its
induced vector field on $M$. The basic forms
on $M$ give rise to a differential complex whose cohomology coincides
with the \v{C}ech cohomology of the topological quotient $M/G$ (see
\cite{k2}).  

\begin{lemma}\labell{L:basic}
Let $G$ be a compact Lie group. 
Let $(N, \ow, \Phi)$ and $(N', \ow', \Phi')$ be neighborhoods of the
zero fibers in six dimensional Hamiltonian
$G$-manifolds such that $0 \in \Phi(N) = \Phi'(N')$.
Assume that $N$ and $N'$ satisfy Condition \ref{E:tech} and 
that their Duistermaat-Heckman functions
coincide.
Then for any $\Phi$-$G$-diffeomorphism $f \colon N \to N'$, 
there exists a
basic one-form $\beta$ such that $d\beta = f^*\ow'-\ow$.
\end{lemma}

\begin{proof}
Consider the closed two-form $\Omega = f^*\ow'-\ow$. Since $f$
commutes with the group action and both $\ow$ and $\ow'$ are
$G$-invariant, $\Omega$ is also invariant. 
Using the fact that $f^*\ow'$ and
$\ow$ have the same moment map $\Phi$ and $\iota(\xi_N)\ow = -d\la
\Phi, \xi\ra$, we have $\iota(\xi_N)\Omega = 0$ for all $\xi \in
\fg$. So $\Omega$ is basic and is a pull-back of a two-form
$\Tilde{\omega}$ on $N/G$.

By Condition (\ref{E:tech}), it suffices to show that
the restriction of $\Tilde{\Omega}$ to the reduced space
$\Phi^{-1}(G\cdot \alpha)/G$ is exact at some regular value $\alpha
\in U$. 
Since
the reduced space is two dimensional, it is enough to show
that the integral of $\Tilde{\Omega}$
over the reduced space is zero. That is, the symplectic volumes of
the reduced spaces are the same. 
This follows from the fact that the Duistermaat-Heckman functions 
coincide at $\alpha$.

So $\Tilde{\Omega} = d \Tilde{\beta}$ for some $\Tilde{\beta}$, 
and we pull back to obtain $\Omega =
\pi^*d\Tilde{\beta} = d(\pi^*\Tilde{\beta}) = d\beta$, where $\beta
= \pi^*\Tilde{\beta}$ is a basic one-form. 
\end{proof}

\comment{check proof of basic form as a pull-back}

\begin{lemma}\labell{L:nondeg}
Let $G = \Su(2)$ or $\So(3)$ act effectively
on a six dimensional manifold $M$.  
Let $\ow_0$ and $\ow_1$ be 
two $G$-invariant symplectic forms on $M$
with the same moment map $\Phi$ and assume that $\ow_0$ and $\ow_1$ 
induce the same orientation. Then the $G$-invariant two-form
$\ow_t = (1-t)\ow_0+t\ow_1$ is nondegenerate for all $0 \leq t \leq 1$
on a neighborhood of the zero fiber $\Phi^{-1}(0)$.
\end{lemma}

\begin{proof}
Nondegeneracy is a local condition. To show that $\ow_t$ is
nondegenerate on a neighborhood of the zero fiber, 
it is enough to show that it is nondegenerate on the zero fiber. 
Consider a point $x \in \Phi^{-1}(0)$ such that
the local model for the orbit $G \cdot x$ is 
$Y = G \times_H (\fh^0 \times \C^n)$. The tangent space at $x$ splits as 
$T_xY = \fg/\fh \times \fh^0 \times \C^n$. The two-form $\ow_t$ at 
$x$ is then of the form 
	\[
	\begin{pmatrix}
	0 & I & 0 \\
	-I & * & * \\
	0 & * & \tilde{\ow_t}
	\end{pmatrix}
	\]
where $\tilde{\ow_t}$ is an $H$-invariant linear symplectic 
form on $\C^n$ and $I$ is the
natural pairing between the tangent space $\fg/\fh$ and its dual
$\fh^0$. So the two-form $\ow_t$ is nondegenerate if and only if the
corresponding $\tilde{\ow_t}$ is nondegenerate.
	\begin{description}
	\item[Case 1] When the local model is $G \times_{\Gamma}
		\R^3$, the two-form $\ow_t = \left(\begin{smallmatrix}
		0 & I \\ -I & * \end{smallmatrix}\right)$ is nondegenerate
		for all $t$.
	\item[Case 2] When the local model is $G \times_H (\R^2
		\times \C)$ or $G \times_H \R^2 \times \C$ with 
		$H = S^1$ or $N(S^1)$, the two-form $\tilde{\ow_t}$ is 
		a linear symplectic form 
		defined on $\C$ and therefore is 
		$A^t dz \wedge d\bar{z}$
		for some constant $A^t$. By definition,  
		$A^t = (1-t) A^0 + tA^1$. Since $\ow_0$ and $\ow_1$
		induce the same orientation, $A^0A^1 > 0$. So
		$A^t$ never vanishes for $0 \leq t \leq 1$, 
		and $\tilde{\ow_t}$ is nondegenerate for $0 \leq t \leq 1$.
	\item[Case 3] When the local model is $\C^3$, 
		we can translate $\tilde{\ow_t}$ 
		to any other point in $\C^3$. 
		In particular we can translate it to an $x$ in
		the zero fiber with $S^1$ stabilizer; see
		Table~\ref{Ta:model}. Case 2 then 
		applies.
	\item[Case 4] The local model $Y$ is $\C^2 \times \C$.
		Let $\ow_0$ and $\ow_1$ be $\Su(2)$-invariant symplectic
		forms on $\C^2 \times \C$ that have the same moment map and induce
		the same orientation. Let $\ow_t = (1-t)\ow_0 + t \ow_1$.
		Then $\ow_t$ can be written as 
		\[
		\begin{split}
		\ow_t &= \sqrt{-1}\sum_{i=1}^3 A_i^t dz_i \wedge
		d\bar{z_i}\\  
		&+ \sqrt{-1}\sum_{i < j} \left(
		B_{ij}^t dz_i\wedge dz_j + 
		\overline{B_{ij}^t}d\bar{z_i}\wedge d \bar{z_j}\right)\\
			&+ \sqrt{-1}\sum_{i < j}\left(
			C^t_{ij}dz_i \wedge d\bar{z_j} +
			\overline{C^t_{ij}}d\bar{z_i}\wedge dz_j \right)
			\end{split}
			\]
		Because $\ow_t$ is invariant under $\Su(2)$ action, it is
			invariant under the transformations 
			$\left(\begin{smallmatrix}
				a & 0 \\
				0 & \bar{a}
			\end{smallmatrix}\right)$,
			$\left(\begin{smallmatrix}
				\cos\theta & \sin\theta \\
				-\sin\theta & \cos\theta 
			\end{smallmatrix}\right)$, and
			$\left(\begin{smallmatrix}
				\cos\theta	& -i\sin\theta  \\
				i\sin\theta & \cos\theta 
		\end{smallmatrix}\right)$ for any $\theta$. It follows that
		$A_1^t = A_2^t$, and $B_{ij}^t = C_{ij}^t = 0$ for $i < j=3$.
			The induced vector fields for the 
			above circle actions are 
			\[
			\begin{aligned}
			X_{\xi_1}
			&= \sqrt{-1}\,(z_1 \frac{\partial}{\pz_1} 
			- \bar{z_1}\frac{\partial}{\partial{\bar{z_1}}}
			- {z_2}\frac{\partial}{\partial{z_2}} 
			+ \bar{z_2}\frac{\partial}{\partial{\bar{z_2}}}) \\
			X_{\xi_2}
			&=  
			 z_2\frac{\partial}{\pz_1}
			+ \bar{z_2}\frac{\partial}{\partial{\bar{z_1}}}
			- z_1 \frac{\partial}{\pz_2} 
			+ \bar{z_1}\frac{\partial}{\partial\bar{z_2}}\\
			X_{\xi_3}
			&= \sqrt{-1}\,(z_2 \frac{\partial}{\pz_1} 
			- \bar{z_2}\frac{\partial}{\partial\bar{z_1}} 
			+ z_1\frac{\partial}{\pz_2}
			- \bar{z_1}\frac{\partial}{\partial{\bar{z_2}}} )
			\end{aligned}
			\]
	respectively. Because $\iota(X_{\xi_i})\ow_t = -d\Phi_Y^{\xi_i}$,
	direct computation shows that 
	the coefficients $A_{1}^t$, $A_{2}^t$, $B_{12}^t$, and $C_{12}^t$ 
	are determined by the moment map and thus independent of $t$. 

	Note $A_{3}^t = (1-t) A_{3}^0 + tA_{3}^1$. Since $\ow_0$ and $\ow_1$
		induce the same orientation, $A_{3}^0A_{3}^1 > 0$. So
		$A_{3}^t$ never vanishes and ${\ow_t}$ is nondegenerate. 	\end{description}
\end{proof}

\begin{proof}[Proof of Proposition \ref{P:symplecto}]
Suppose $f$ is a $\Phi$-$G$-diffeomorphism from $(M, \ow, \Phi)$ to 
$(M', \ow', \Phi')$. By Lemma
\ref{L:basic}, there exists a basic one-form $\beta$ such that
$d\beta = f^*\ow'-\ow$. Consider $\ow_t = (1-t)\ow + t f^*\ow'$ for $0
\leq t \leq 1$. Then $\ow_0 = \ow$ and $\ow_1 = f^*\ow'$ have the same
moment map $\Phi$ and induce the same orientation. By Lemma \ref{L:nondeg},
there exist a neighborhood of the zero fiber 
$N \subset M$ on which $\ow_t$ is nondegenerate. 
We can solve $\iota(X_t)\ow_t =
-\beta$ for the time dependent vector field $X_t$. Denote $\varphi^t$
the flow of this vector field $X_t$ satisfying $\varphi^0 =
\text{id}$. Because $\la d\Phi(X_t), \xi \ra = -\ow_t(\xi_M, X_t) =
(\iota(X_t)\ow_t)(\xi_M) = -\beta(\xi_M) = 0$, the flow $\varphi^t$
preserves the fibers of the moment map $\Phi$. 
Since $\Phi$ is proper,
$\varphi^t$ exists for all $0 \leq t \leq 1$. 

Define $F_t = f \circ \varphi^t$. Since $F_t^*\Phi' = \Phi'\circ F_t 
= \Phi' \circ f \circ \varphi^t = \Phi \circ \varphi^t = \Phi$, we
know that $F_t\colon N\to N'$ respects the moment maps for $N' = f(N)$. 
Since $\ow_t$ and $\beta$ are invariant,
$X_t$ is invariant. So $\varphi^t$ and hence $F_t$ is $G$-equivariant.  
Finally
\[
\begin{split}
\frac{d}{dt}\left( \varphi^* \ow_t \right) 
&=  \varphi^*_t
L_{X_t}\ow_t + \varphi^*_t \left( \frac{d}{dt}\ow_t\right) \\
&= \varphi^*_t \left( L_{X_t}\ow_t+\frac{d}{dt}\ow_t\right) \\
&= \varphi^*_t \left( d\iota(X_t)\ow_t -\ow +f^*\ow' \right)\\
&= \varphi^*_t \left( d(-\beta) +d\beta \right) \\
&=0 
\end{split}
\]
So $\varphi^*_1\ow_1 = \varphi^*_0 \ow_0 = \varphi^*_1 f^*\ow' = (f
\circ \varphi_1)^* \ow' = F_1^*\ow'= \ow$ and $F_1\colon  N \to N'$ is 
an equivariant symplectomorphism that respects the moment maps.  
\end{proof}

\section{Passing to the quotient}\labell{S:toquotient}

In this section, we show that it is enough to work with
a specific class of diffeomorphisms of the quotients rather 
than $\Phi$-$G$-diffeomorphisms. 

Let $G=\Su(2)$ or $\So(3)$ act on a manifold $M$. With quotient
topology, the quotient $M/G$ has a natural smooth structure; a
function is smooth if and only if its pull-back to $M$ is smooth. 
A \textbf{smooth} function between quotients is a map $f\colon M/G \to
M'/G$ that pulls back smooth functions to smooth functions. A smooth
function $f$ between quotients 
is a \textbf{diffeomorphism} if it is smooth and has a 
smooth inverse. If $M$ and $M'$ are oriented, the smooth parts of the
quotients $M/G$ and $M'/G$ can be oriented with a choice of the 
orientation on $G$. Whether or not a diffeomorphism $f\colon M/G \to M'/G$
preserves the orientation is independent of that choice. 

While this notion of diffeomorphism is natural, 
we will use a stronger notion of $\Phi$-diffeomorphisms which allows
us to have a better control over neighborhoods of the exceptional orbits. 

\comment{Can this be defined for arbitrary $G$? Check Kirwan}

First we observe that when $G = \Su(2)$ or $\So(3)$, we can identify
the dual of its Lie algebra $\fg^*$ with $\R^3$ and its orbit space
$\fg^*/G$ under
the coadjoint action with $\R_+ \simeq [0, \infty)$. 
Since the moment map $\Phi$ is equivariant, it induces a map
$\overline{\Phi}\colon M/G \to \R_+$ such that the following diagram commutes:
\begin{equation}\labell{E:norm2}
\begin{diagram}
\node{M} \arrow{e,t}{\Phi}\arrow{s,l}{\pi}\node{\fg^* = \R^3}\arrow{s,r}{p}\\
\node{M/G} \arrow{e,t}{\overline{\Phi}}\node{\fg^*/G = \R_+}
\end{diagram}
\end{equation}
where $p$ is the norm square $|\xi|^2$ and
$\overline{\Phi}([m])=|\Phi(m)|^2$.  

\comment{need equivalent definition for the cross-section}

\begin{definition}\labell{D:phidiffeo}
Let $G$ be $\Su(2)$ or $\So(3)$.
Let $M$ and $M'$ be oriented manifolds with $G$ actions and
$G$-equivariant maps $\Phi\colon M \to \fg^*$ and 
$\Phi'\colon  M' \to \fg^*$. 
%Let $\overline{\Phi}$ and
%$\overline{\Phi'}$ be induced by the norm squares of the moment maps. 
A \textbf{$\mathbf{\Phi}$-diffeomorphism} 
from $M/G$ to $M'/G$ is an orientation
preserving diffeomorphism $\psi\colon  M/G \to M'/G$ such that 
\begin{enumerate}
\item $\psi^* \overline{\Phi'} = \overline{\Phi}$.
\item $\psi$ and $\psi^{-1}$ lift to a 
$\Phi$-$G$-diffeomorphism in a
neighborhood of each exceptional orbit. 
\end{enumerate}
\end{definition}

We start with a series of lemmas. 

\begin{lemma}\labell{L:locallift}
Let $G$ be $\Su(2)$ or $\So(3)$ and $(M, \ow, \Phi)$ be a 
six dimensional Hamiltonian $G$-manifold.  
Let $Y$ be a local model for a nonexceptional orbit in the zero fiber
of the moment map. Let $W$ and $W'$ be invariant
open subsets of $Y$. Let $f\colon W/G \to W'/G$ be a $\Phi$-diffeomorphism.
Then $f$ lifts to a $\Phi$-$G$-diffeomorphism $F\colon W \to W'$. 
\end{lemma}

\begin{proof}
When $G = \Su(2)$, the only possible local model for a 
nonexceptional orbit in the zero fiber is 
$Y = \C^2 \times \C$. The quotient $Y/G$ can be identified as $\C^2/G
\times \C$. The $\Phi$-diffeomorphism $f\colon Y/G \to Y'/G$ necessarily
has the form 
$f([w], z) = ([w], \varphi([w], z))$ for some diffeomorphism
$\varphi\colon  \C^2/G \times \C \to \C$. 
It locally lifts to a $\Phi$-$G$-diffeomorphism $F(w, z)
= (w, \varphi([w], z))$ for $w \in \C^2$ and $z\in \C$.  

When $G = \So(3)$, the local model for a nonexceptional orbit is 
$G \times_H \R^2\times \C$ with $H = S^1$ or $N_G(S^1)$. The
quotient $Y/G$ can be identified as $\R^2/H \times \C$. Then the
$\Phi$-diffeomorphism $f\colon  Y/G \to Y'/G$ can be written as
 $f([\mu], z) = ([\mu],
\varphi([\mu], z))$ for some $\varphi\colon  \R^2/H \times \C \to \C$. 
It locally lifts to a $\Phi$-$G$-diffeomorphism $F([g, \mu], z)
= ([g, \mu], \varphi([\mu], z))$.  
\end{proof}

\begin{lemma}\labell{L:unitaryfij}
Let $G = \Su(2)$ act effectively on a six dimensional symplectic 
manifold $(M, \ow)$. Assume the action is Hamiltonian and the moment
map is $\Phi$. 
Consider the local model $Y = \C^2 \times \C$ for 
a nonexceptional orbit in the zero 
fiber. 
Let $F\colon Y \to Y$ be an equivariant diffeomorphism 
that preserves the orbits and respects the moment maps. 
Extend the $\Su(2)$ action to $\operatorname{U}(2)$.
Then there exists a smooth invariant 
function $h\colon Y \to S^1$ such that $F(y) = h(y) \cdot y$ where the $S^1$
action on $Y$ is induced by $\operatorname{U}(2)$ on $\C^2$.
\end{lemma} 

\comment{better way to define the action?}

\begin{proof}
Let $F\colon Y \to Y$ be an equivariant diffeomorphism that preserves the
orbits and respects the moment maps. 
Identify the local model for a nonexceptional orbit as $Y = \C^2
\times \C$.
Then $F$ takes $(w, z)$ to some $(w', z')$ and  
the following diagram commutes,
\[
\begin{diagram}
\node{\C^2 \times \C} \arrow{e,t}{F} \arrow{s,l}{\pi}
\node{\C^2 \times \C}\arrow{s,r}{\pi}\\
\node{\C^2 / G \times \C} \arrow{e,t}{\text{id}}\node{\C^2 /G \times \C}
\end{diagram}
\]
where $\pi(w, z) = ([w], z)$. Therefore $z' = z$, and $w' = hw$ for
some $h \in \operatorname{U}(2)$. Since $F$ is equivariant,
$hgw=ghw$ for all $g$ in $\Su(2)$. And hence $h$ belongs to the
center of $\operatorname{U}(2)$, i.e., $ h \in S^1$.
\end{proof}

\comment{this proof is not clear}

\begin{lemma}\labell{L:fixedsurface}
Let $G$ denote $\Su(2)$.
Let $(N, \ow, \Phi)$ and $(N', \ow', \Phi')$ be 
neighborhoods of tall zero fibers 
in six dimensional Hamiltonian $G$-manifolds.  
Assume that $N$ and $N'$ satisfy Condition
\ref{E:tech} and their
Duistermaat-Heckman functions are the same. 
Then every $\Phi$-diffeomorphism from
$N/G \to N'/G$ that locally lifts to a $\Phi$-$G$-diffeomorphism 
globally lifts to a $\Phi$-$G$-diffeomorphism.
\end{lemma}

\comment{If we use local models in the proof, in general it is not
enough to assume that the zero neighborhood $N$ satisfies Condition
\ref{E:tech}. We need to shrink $N$ so that we can cover it by local
models.} 

\begin{proof}
Assume ${\psi}\colon N/G \to N'/G$ is a $\Phi$-diffeomorphism. 
Choose an open invariant cover $\mathcal{U}$ of $N$ such that 
$U_i \cap \Phi^{-1}(0) \neq \emptyset$ for each $U_i \in
\mathcal{U}$. Take a refinement if necessary, we can assume that
each $U_i$ is an invariant open subset of the local model 
$\C^2 \times \C$; see Table \ref{Ta:model}. 
By Lemma \ref{L:locallift}, 
$\psi$ locally lifts to a $\Phi$-$G$-diffeomorphism 
$\Psi_i\colon U_i \to N'$ for each $U_i$. 
By Lemma \ref{L:unitaryfij}, there exists smooth invariant functions 
$f_{ij}\colon  U_i \cap U_j \to S^1$ such that $\Psi_i = f_{ij} \cdot
\Psi_j$. The collection of $\{f_{ij}\}$ form a \v{C}ech cocycle $g$ on
$\mathcal{U}$ with coefficient in $S^1$. If the corresponding class
$[g] \in \check{H}^1(N/G, S^1)$ is trivial, ${\psi}$ lifts
to a global $\Phi$-$G$-diffeomorphism. 

The short exact sequence $0 \to \Z \to \R \to S^1 \to 0$ determines a
long exact sequence in cohomology. Since there exists a smooth
partition of unity on $N/G$, $\check{H}^i(N/G; \R) = 0$ for all $i >
0$. In other words, $\check{H}^1(N/G; S^1) = \check{H}^2(N/G; \Z)$. 
According to Condition \ref{E:tech}, the restriction 
$\check{H}^2(N/G; \Z) \to \check{H}^2(\Sigma; \Z)$ is one-to-one, 
where $\Sigma =\Phi^{-1}(G \cdot \alpha)/{G} = \Phi^{-1}(\alpha)/{S^1}$ 
is a regular reduced space; we
only need to show that the image of $[g]$ in $\check{H}^2(\Sigma; \Z)$ 
vanishes. 
Let $\{\lambda_j\}$ be a partition of unity subordinate to the cover
$\mathcal{U} \cap \Phi^{-1}(\alpha)$.  
The \v{C}ech-de Rham isomorphism takes the image of $[g]$ to
the cohomology class of the basic differential two-form 
\begin{equation}\labell{E:chernclass}
\Sigma_j d\lambda_j g_{ij}^{-1}dg_{ij}
\end{equation}
when restricted to the open set $U_i \cap \Phi^{-1}(\alpha)$.
It is equal to the difference
between the curvature form $d\Theta$ and the pullback
$\Psi^*d\Theta'$, where
$\Theta'$ is a connection one-form on ${\Phi'}^{-1}(\alpha)$ and
$\Theta = \sum \lambda_i \Psi^*_i\Theta'$ is a connection one-form on
$\Phi^{-1}(\alpha)$. The integrals of the curvature forms over the
reduced spaces are equal to the slopes of the
Duistermaat-Heckman functions at $\alpha$. 
Since the Duistermaat-Heckman functions
are the same, (\ref{E:chernclass}) is exact as a basic form and $[g]=0$.
\end{proof}

\comment{explain better}

\begin{lemma}\labell{L:fij}
Let $G = \So(3)$ act effectively on a six dimensional symplectic 
manifold $(M, \ow)$. Assume the action is Hamiltonian and the moment
map is $\Phi$. 
Consider the local models for nonexceptional orbits of the zero fiber 
$Y = \So(3)\times_H \R^2 \times \C$ with $H = S^1$ or $N(S^1)$. 
Let $F\colon Y \to Y$ be an equivariant diffeomorphism 
that preserves the orbits and respects the moment maps. 
Then there exists a smooth invariant 
function $h\colon Y \to N(S^1)/H$ such that $F(y) = h(y) \cdot y$ where the
action is induced by the extension of the $H$ action on $\So(3)
\times \R^2 \times \C$ to an action of $N(S^1)$.
\end{lemma}

\begin{proof}
Let $F\colon Y \to Y$ be an equivariant diffeomorphism that preserves the
orbits and respects the moment maps. 
Identify the local model for a nonexceptional orbit as $Y = (\So(3)
\times_H \R^2) \times \C$ where $H =S^1$ or $N(S^1)$. 
Then $F$ takes $([g, \mu], z)$ to some $([g', \mu'], z')$ and  
the following diagram commutes,
\[
\begin{diagram}
\node{(\So(3) \times_H \R^2) \times \C} \arrow{e,t}{F} \arrow{s,l}{\pi}
\node{(\So(3) \times_H \R^2) \times \C}\arrow{s,r}{\pi}\\
\node{\R^2 / H \times \C} \arrow{e,t}{\text{id}}\node{\R^2/H \times \C}
\end{diagram}
\]
where $\pi([g, \mu], z) = ([\mu], z)$. Since $F$ is equivariant,
$g'=gh$ for some $h$ in $\So(3)$. $F$ respects the moment maps, so
$Ad^{\dag}(g)\mu = Ad^{\dag}(g')\mu'$ and $\mu' = {g'}^{-1}g\mu =
h^{-1}g^{-1}g\mu=h^{-1}\mu$. Since $[\mu'] = [\mu] = [h^{-1}\mu]$, 
$h \in N(S^1)$. Now $[gh,
h^{-1}\mu] = [g, \mu]$ if $h \in H$. 
Hence $F([g, \mu], z) = ([gh, h^{-1}\mu], z)$ where $h\colon  Y \to N(S^1)/H$.
\end{proof}

\begin{definition}\labell{D:orientbundle}
For any real $n$ dimensional vector bundle $\pi\colon W \to X$, 
there exists an associated orientation bundle
$p\colon \Tilde{X} \to X$, whose fiber over a point $x$ is the 
two ways to orient $\pi^{-1}(x)$. This is a two-sheeted covering, and 
\v{C}ech cocycles provide a convenient way to construct it. 
Choose an open cover $\mathcal{U} = \{U_\alpha\}$ of $X$ with
trivialization maps 
$\varphi_\alpha\colon  U_\alpha \times \R^n \to \pi^{-1}(U_\alpha)$. 
The Jacobian determinants of the change of fiber coordinates from 
$\R^n$ to $\R^n$ have a locally constant sign, 
which gives a locally constant function from 
$U_\alpha \cap U_\beta$ to $\Z_2$. The chain rule for Jacobians
implies that this is a cocycle. 
%Define $p \colon  \Tilde{X} \to X$ to be
%the resulting $\Z_2$ bundle. This is the orientation bundle
%associated with the vector bundle $W$ and 
It determines
an element $w_1(W) \in H^1(X;\Z_2)$, called the \textbf{first
Stiefel-Whitney class} of $W$. In a similar fashion, 
we can construct the associated orientation bundle and 
define the first Stiefel-Whitney class of any fiber bundle $W$ 
over a manifold $X$ provided that the fiber of $W$ is 
connected and orientable. 

In particular, when the zero fiber of a Hamiltonian $\So(3)$-manifold
is tall, and when the principal stabilizer of the zero fiber 
is $S^1$, the zero fiber $\Phi^{-1}(0)$ is
a sphere bundle over the reduced space 
off the exceptional orbits. Let $E = \{\,E_j\,\}$ denote the set of
exceptional orbits in the zero fiber.
Let $M^{\text{reg}}_0$ denote the smooth part of the symplectic
quotient at $0$, i.e., $M_0^{\text{reg}}=(\Phi^{-1}(0) \smallsetminus
E)/G$. Then $M_0^{\text{reg}}$ is diffeomorphic to 
$\Sigma \smallsetminus \{\text{\emph{finitely many points}}\}$,
where $\Sigma = \Phi^{-1}(0)/G$ is a closed connected oriented
surface.  
Through the orientations on the fiber spheres,
$\Phi^{-1}(0)$ 
induces an associated orientation bundle on $M_0^{\text{reg}}$,
and the first Stiefel-Whitney class in $H^1(M_0^{\text{reg}};\Z_2)$.
\end{definition}

\begin{lemma}\labell{L:nonmixed}
 Let $G$ denote $\So(3)$. 
Let $(N, \ow, \Phi)$ and $(N', \ow', \Phi')$ 
be neighborhoods of zero fibers which satisfy Condition \ref{E:tech} 
in six dimensional Hamiltonian $G$-manifolds.
Assume that the zero fibers are tall with the same principal
stabilizer, and that they have no exceptional orbits. 
Then every $\Phi$-diffeomorphism from
$N/G \to N'/G$ that locally lifts to a $\Phi$-$G$-diffeomorphism 
globally lifts to a $\Phi$-$G$-diffeomorphism if one of the following
condition holds:
\begin{enumerate}
\item the principal stabilizer of the zero fibers is $N(S^1)$; or
\item the principal stabilizer of the zero fibers is $S^1$ and 
the first Stiefel-Whitney class on $\Phi^{-1}(0) \to
{\Phi^{-1}(0)}/{G}$ equals the pull-back of the 
first Stiefel-Whitney class on 
${\Phi'}^{-1}(0) \to {\Phi'}^{-1}(0)/{G}$.
\end{enumerate}
\end{lemma}

\begin{proof}
Assume $\psi\colon N/G \to N'/G$ is a $\Phi$-diffeomorphism. 
Choose an open invariant cover $\mathcal{U}$ of $N$ such that $U_i
\cap \Phi^{-1}(0) \neq \emptyset$ for each $U_i \in \mathcal{U}$. 
Take a refinement if necessary, we can assume that each $U_i$ is an
invariant open subset of the local model $G\times_H \R^2 \times \C$
where $H$ is $S^1$ or $N(S^1)$; see Table \ref{Ta:model}. 
By Lemma \ref{L:locallift},  
$\psi$ locally lifts to a $\Phi$-$G$-diffeomorphism 
$\Psi_i\colon U_i \to N'$. 
By Lemma \ref{L:fij}, there exists smooth invariant functions $f_{ij}\colon 
U_i \cap U_j \to N(S^1)/H$ such that $\Psi_i = f_{ij} \cdot \Psi_j$ where
$H = S^1$ or $N(S^1)$. 
     
\begin{description}
\item[Case 1] $H = N(S^1)$. Then $f_{ij}$ is the 
identity for all $i, j$ 
and $\Psi_i = \Psi_j$ on $U_i \cap U_j$. 
Therefore $\psi$ lifts globally. 
\item[Case 2] $H = S^1$. The set $\{f_{ij}\}$ defines
a \v{C}ech cocycle $f \in
\Check{C}^1(\mathcal{U}; \Z_2)$. 
Then there exists a global lift if the
class $[f] \in \Check{H}^1(N/G; \Z_2)$ is trivial. 
By Condition \ref{E:tech}, the restriction 
$H^1(N/G; \Z_2) \to
H^1(\Phi^{-1}(0)/G; \Z_2)$ is one-to-one. It suffices to show that
the image of $[f]$ in $H^1(\Phi^{-1}(0)/G; \Z_2)$ is trivial. 

Let $X = \Phi^{-1}(0)/G$ and $X' = {\Phi'}^{-1}(0)/G$.
Let $\pi\colon  \Tilde{X} \to X$ and 
$\pi'\colon  \Tilde{X'} \to X'$ be the associated orientation
bundles described in Definition \ref{D:orientbundle}. Refine the cover if
necessary so that we have a good cover (simply connected, locally
path connected, and each intersection is connected), again denoted 
by $\mathcal{U}$,
on $\Phi^{-1}(0)/G$ and $\mathcal{U'} = \psi(\mathcal{U})$ 
also a good cover on ${\Phi'}^{-1}(0)/G$. 
Let $\varphi_i\colon  \pi^{-1}(U_i) \to U_i
\times \Z_2$ and $\varphi'_i\colon  {\pi'}^{-1}(U_i') \to U_i' \times \Z_2$ be 
the trivialization maps for these two bundles subject to the good covers. 
Let $\{g_{ij}\}$ and $\{g_{ij}'\}$ be the transition functions for
$\{\varphi_i\}$ and $\{\varphi'_i\}$. Then $\{g_{ij}\}$ is a \v{C}ech cocycle for 
the first Stiefel-Whitney class 
$w_1(\Phi^{-1}(0))$ and $\{g_{ij}'\}$ for $w_1({\Phi'}^{-1}(0))$. 

The proof of Lemma \ref{L:fij} implies that there exists a 
constant function $h_i \colon  U_i \to \Z_2$ induced by the 
restriction of a local $\Phi$-$G$-diffeomorphism $\Psi_i$ to the zero 
fiber. $h_i$ maps $(x, t) \in U_i \times \Z_2$, a local chart of the
associated orientation bundle, to 
$(\psi(x), h_i + t) \in U_i' \times \Z_2$. Then 
$h_i + g_{ij}'\psi - h_j - g_{ij} = f_{ij}$. The class $[f]$
vanishes exactly when the first 
Stiefel-Whitney classes are the same.
\end{description}
\end{proof}

\comment{rewrite, can we use $\Z$ coefficient in Condition
\ref{E:tech} and use the universal coefficient theorem
without assuming torsion free?}

\begin{example} Consider any Riemann surface $\Sigma$ of genus $k$,
and any $w_1 \in H^1(\Sigma; \Z_2)$. Let $P \to \Sigma$
be a principal $\Z_2$ bundle whose first Stiefel-Whitney class
is $w_1$. Consider the associated bundle 
$M = P \times_{\Z_2} T^*S^2$ where 
the $\Z_2$ action on $T^*S^2$ is the lifted  
action induced by the antipodal map on the sphere. 
There exists a symplectic form on $M$ and a corresponding moment map.
Up to equivariant diffeomorphisms, $M$ is
determined by $w_1$ and $k$. The zero fiber is $P \times_{\Z_2} S^2$,
a sphere bundle that determines $w_1$, and the quotient $M/\So(3) =
\Sigma \times \R_+$ determines the genus. 
\end{example}

%\begin{remark}
%If every stabilizer of the zero fiber is $S^1$, 
%the local model $\So(3) \times_{S^1} \R^2 \times \C$ implies that 
%$\Phi^{-1}(0)$ is locally $S^2 \times \C$. Therefore the zero
%fiber is topologically a $S^2$ bundle. The action of $\So(3)$ maps
%each fiber $S^2$ into itself, and this action is transitive. Then the
%transition maps of this sphere bundle can be chosen so that they take
%their values in $N(S^1)/S^1 \simeq \Z_2$ \cite{ls}. 
%Consequently there exists a
%principal $\Z_2$ bundle $P$ over $\Phi^{-1}(0)/\So(3)$ and
%$\Phi^{-1}(0) \simeq P \times_{\Z_2} S^2$, where $Z_2$ acts on
%$S^2$ by antipodal maps and the neighborhood of $\Phi^{-1}(0)$ is 
%$P \times_{\Z_2} (T^*S^2)$.     
%\end{remark}

\begin{lemma} \labell{L:mixed}
Let $G$ denote $\So(3)$. 
Let $(N, \ow, \Phi)$ and $(N', \ow', \Phi')$ be 
neighborhoods of zero fibers which satisfy Condition \ref{E:tech} 
in six dimensional Hamiltonian $G$-manifolds. Assume that the zero
fibers are tall and that there exist exceptional
orbits $E= \{E_j\}$ on the zero fiber of $\Phi$ and $E'=\{E_j'\}$ of 
$\Phi'$. Every $\Phi$-diffeomorphism from
$N/G \to N'/G$ that locally lifts to a $\Phi$-$G$-diffeomorphism 
globally lifts to a $\Phi$-$G$-diffeomorphism 
if the first Stiefel-Whitney class for
$\Phi^{-1}(0) \smallsetminus E 
\to \frac{\Phi^{-1}(0) \smallsetminus E}{G}$ 
equals the pull-back of the first Stiefel-Whitney class for
${\Phi'}^{-1}(0) \smallsetminus E' \to \frac{{\Phi'}^{-1}(0)
\smallsetminus E'}{G}$.
\end{lemma}

\begin{proof}
The assumption that there exist exceptional orbits implies that
the principal stabilizers of the zero fibers are $S^1$. After 
passing to a suitable finer covering $\mathcal{U} = \{U_i\}$ 
if necessary, we can assume that on $\Phi^{-1}(0)/G$, each 
$U_i$ is simply connected and locally path connected, 
and that $U_i \cap U_j$ contains no exceptional orbits. We 
then apply a similar argument as before on 
$\Phi^{-1}(0) \smallsetminus E$. 
\end{proof}

\begin{remark}
There are two types of exceptional orbits that can occur in a tall
zero fiber with $S^1$ stabilizer.

The first type of the exceptional orbits are the isolated fixed 
points in $\Phi^{-1}(0)$. It corresponds to
$\{0\}$ in the local model $Y = \C^3 \simeq T^*\R^3 = \{(q, p)\}$. 
Define the map 
\[
\eta\colon  \Phi_Y^{-1}(0) \smallsetminus \{0\} \to 
\frac{\Phi_Y^{-1}(0) \smallsetminus \{0\}}{\So(3)} \to 
\frac{\R^2 \smallsetminus \{0\}}{\Z_2} 
\to \R^2 \smallsetminus \{0\}
\]
by $(q, p) \mapsto \text{orbit}(q, p)=\text{orbit}(u,0,0,v,0,0)
\mapsto [u,v]=[r,\theta] \mapsto (r, 2\theta)$. Let $I$ denote the
interval $[0, \pi]$. Then there exists a commutative diagram 
\[
\begin{diagram}
\node{I\times S^2}\arrow{s,l}{\text{pr}}\arrow{se,t}{f}\\
\node{\frac{I\times S^2}{(0,x) \sim (\pi,
-x)}}\arrow{e,b}{g}\node{\eta^{-1}(S^1)}
\end{diagram}
\] where $f(t, x) = ((\cos t)x, (\sin t)x)$ induces a homeomorphism $g$. 
In other words, the product of the values of the 
$\{f_{ij}\}$ defined in Lemma \ref{L:fij} as one moves along a loop
around $E/G$ in $\Phi^{-1}(0)/G$
is the nontrivial element in $\Z_2$. 

The second type of the exceptional orbits are those whose stabilizer 
is $N(S^1)$ and whose local model is $\So(3) \times_{N(S^1)}(\R^2 \times \C)$.
The local model immediately
implies that the product of the values of the 
$\{f_{ij}\}$ defined in Lemma \ref{L:fij} as one moves along a loop
around $E/G$ in $\Phi^{-1}(0)/G$
is the nontrivial element in $\Z_2$. 

In fact, blowing up an isolated fixed point gives us an exceptional orbit
of the second type. 
\end{remark}

\begin{proposition}\labell{P:globallift} 
Let $G$ denote $\Su(2)$ or $\So(3)$. Let 
$(N, \ow, \Phi)$ and $(N', \ow', \Phi')$ be neighborhoods of tall
zero fibers in six dimensional Hamiltonian $G$-manifolds. Assume that
$N$ and $N'$ satisfy Condition \ref{E:tech} and that their 
Duistermaat-Heckman functions are the same. Assume that
the zero fibers $\Phi^{-1}(0)$ and ${\Phi'}^{-1}(0)$ have the same
principal stabilizer and the same 
first Stiefel-Whitney
class when applicable. Then   
every $\Phi$-diffeomorphism from $N/G$ to $N'/G$ globally lifts 
to a $\Phi$-$G$-diffeomorphism. 
\end{proposition}

\begin{proof}
This is a direct result from 
Lemmas \ref{L:locallift}, \ref{L:fixedsurface}, \ref{L:nonmixed}, and
\ref{L:mixed}.  
\end{proof}

\section{The topology of the quotient} \labell{S:qtop}

Let $G$ be $\Su(2)$ or $\So(3)$. Let 
$(N, \ow, \Phi)$ and $(N', \ow', \Phi')$ be neighborhoods of tall
zero fibers that satisfy Condition \ref{E:tech} 
in six dimensional Hamiltonian $G$-manifolds.
By Propositions \ref{P:symplecto} and \ref{P:globallift}, we have shown
that as long as $N$ and $N'$ have the same Duistermaat-Heckman function, 
and their zero fibers have the same principal stabilizer as well as the
same first Stiefel-Whitney class,
$N$ and $N'$ are isomorphic if their quotients $N/G$ and $N'/G$ are 
$\Phi$-diffeomorphic.  In the following sections, we will show that
the isotropy data and the genus
determine the quotients up to  
$\Phi$-diffeomorphisms. First, we describe the topology of the quotients.

\begin{proposition}\labell{P:qtop}
Let $G$ be $\Su(2)$ or $\So(3)$ and let $(N, \ow, \Phi)$ be a
$G$-invariant neighborhood of the zero fiber in a six
dimensional Hamiltonian $G$-manifold.
Assume that the zero fiber is tall. 
The quotient $N/G$ is topologically a manifold with boundary. 
In particular, each reduced space $\Phi^{-1}(G \cdot \alpha)/G$ 
is topologically a closed connected oriented surfaces for
$\alpha \in \Phi(N)$. 
\end{proposition}

\begin{proof}
By the Local Normal Form Theorem \ref{T:normalform} and 
Lemma \ref{L:homeo} below, the
quotient $N/G$ is topologically a manifold with boundary. 

By the Local Normal Form Theorem \ref{T:normalform} and 
Corollary \ref{C:alpha-homeo}, the reduced space is a
topological surface. Since the moment map is proper and 
since every moment fiber is connected, 
$\Phi^{-1}(G\cdot\alpha)/G$ is closed and
connected. It is oriented 
since the symplectic structure induces an orientation. 
\end{proof}

\begin{lemma} \labell{L:homeo}
Let $G$ be $\Su(2)$ or $\So(3)$. Identify $\fg^*/G$ with $\R_+$.
Let $(M, \ow, \Phi)$ be a six dimensional Hamiltonian $G$-manifold
with a tall zero fiber. For any local model $Y$ of an orbit in
$\Phi^{-1}(0)$, 
there exists a map 
\[
F_Y=(\overline{\Phi}_Y, \overline{P}_Y)\colon  Y/G \to \R_+ \times \C
\] 
homeomorphic into its image, 
where $\overline{\Phi}_Y\colon Y/G \to \R_+$ is 
induced by the 
moment map, and $\overline{P}_Y$ 
is induced from a $G$-invariant map $P_Y\colon  Y \to \C$. 
\end{lemma}

We call such an $F_Y$ a \textbf{trivializing homeomorphism} of the local
model $Y$. 
%$P_Y$ is called the \textbf{defining polynomial} of the local model $Y$. 

\begin{corollary} \labell{C:alpha-homeo}
The restriction of $F_Y$ to the reduced space is a homeomorphism
for all $\alpha$ such that $|\alpha|^2 \in
\text{image}\,\overline{\Phi}_Y = \text{image}\,|\Phi_Y|^2$. 
\end{corollary}

\begin{proof}[Proof of Lemma \ref{L:homeo}] 
Table \ref{Ta:trivial} gives one collection of trivializing
homeomorphisms for all possible local models when the zero fiber is
tall (cf.~Table \ref{Ta:model}). It is easy to verify that
each $F_Y$ consists of the norm 
square of the moment map and a map $\overline{P}_Y\colon Y/G \to \C$ 
induced from a $G$-invariant map $P_Y\colon  Y \to \C$.  
\begin{table}[h]
\centering
\begin{tabular}{|l|l|}\hline
\emph{Local model $Y$} & 
\emph{Trivializing homeomorphism $F_Y$} \\ \hline
$\So(3)\times_{S^1} \R^2 \times \C$ & 
$F_Y[[g, \mu], z] = (|\mu|^2, z)$ \\ \hline
$\So(3) \times_{N(S^1)}\times \R^2 \times \C$ & 
$F_Y[[g, \mu], z] = (|\mu|^2, z)$ \\ \hline 
$\So(3)\times_{N(S^1)} (\R^2 \times \C)$ & 
$F_Y[[g, \mu, z]] = (|\mu|^2, z^2)$ \\ \hline
$\C^3 = \{\, q+\sqrt{-1}\,p\,\}$ & 
$F_Y[q, p] = (|q \times p|^2, |q|^2 - |p|^2, 2\la q, p\ra)$ \\ \hline 
$\C^2 \times \C$ & 
$F_Y[w, z] = (\frac{1}{4}|w|^2, z)$ \\ \hline
\end{tabular}
\vspace{.2cm}
\caption{Trivializing homeomorphisms.}\labell{Ta:trivial}
\end{table}

For the case 
$Y = \C^3 = \{\,(x, y, z) = q + p\sqrt{-1}\,|\,x, y, z \in \C, \, q =
\text{Re}(x, y, z),\, p=\text{Im}(x, y, z)\,\}$, 
it is easier to consider $F_Y$ as a composition of two functions.
Consider the map $F_1\colon \C^3 / \So(3) \to \R^3$ given by
$F_1([q, p]) = (|q|^2, |p|^2, \la q, p \ra)$
and the map $F_2\colon  \R^3 \to \R_+ \times \C$ given by 
$F_2(\alpha, \beta, \gamma) =
(\alpha\beta - \gamma^2, \alpha-\beta, 2\gamma)$. 
We define $F_Y = F_2 \circ F_1$ so that 
$F_Y([q, p]) = (|q\times p|^2,\, |q|^2-|p|^2, 2\la q, p \ra)$, 
or equivalently, $F_Y([x, y, z]) = 
(|\la y\bar{z}-\bar{y}z, z\bar{x}-\bar{z}x, 
x\bar{y}-\bar{x}y \ra|^2,\, x^2+y^2+z^2)$.

The component functions of $F_1$
generate the $\So(3)$-invariant functions on $\C^3$. 
It is one-to-one, continuous and proper. Hence it is a homeomorphism 
from $\C^3/\So(3)$ into its image in $\R^3$, which is a
solid half cone 
$\alpha \geq 0, \beta \geq 0, \alpha\beta \geq \gamma^2$. 
Since $F_2^{-1}(a, b, c) = 
(\frac{1}{2}(b+\sqrt{4a+b^2+c^2}), \frac{1}{2}(-b+\sqrt{4a+b^2+c^2},
\frac{c}{2})$ is continuous, 
$F_2$ is a homeomorphism from the solid cone to its
image in $\R_+ \times \C$. 
Therefore $F = F_2 \circ F_1$ is also a homeomorphism. 

In each of the other cases, routine checks show that $F_Y$
is well-defined, bijective, continuous and proper. It follows  
that $F_Y$ is a homeomorphism.
\end{proof}

\section{The smooth structure on the quotient}\labell{S:qsmooth}

In this section we study the smooth structure on the quotient. By a
theorem of Schwarz \cite{sCh1}, any invariant smooth function can be
expressed as a smooth function of real invariant polynomials. Using
this fact, we show that the trivializing homeomorphism defined in
Table \ref{Ta:trivial} in the previous section  
is a diffeomorphism on the complement of the exceptional orbits.

\comment{rephrase Schwarz}

Let $G$ be $\Su(2)$ or $\So(3)$ and let $(M, \ow, \Phi)$ be a six
dimensional Hamiltonian $G$-manifold with a tall zero fiber. 
First we list all the exceptional orbits in the local models of an
orbit in the zero fiber: 

\begin{table}[h]
\centering
\begin{tabular}{|l|l|}\hline
\emph{Local model} & 
\emph{Exceptional orbits} \\ \hline
$\So(3)\times_{S^1} \R^2 \times \C$ & 
none \\ \hline
$\So(3) \times_{N(S^1)}\times \R^2 \times \C$ & 
none \\ \hline 
$\So(3)\times_{N(S^1)} (\R^2 \times \C)$ & 
$\{\,[g, \mu, 0]\,\}$ \\ \hline
$\C^3$ & 
$\{0\}$ \\ \hline 
$\C^2 \times \C$ & 
none \\ \hline
\end{tabular}
\vspace{.2cm}
\caption{Exceptional Orbits.}\labell{Ta:except}
\end{table}

\begin{lemma} \labell{L:diffeo}
Let $G$ be $\Su(2)$ or $\So(3)$ and $(M, \ow, \Phi)$ be a six
dimensional Hamiltonian $G$-manifold such that $0 \in \Phi(M)$. 
Assume that the zero fiber is tall. 
Let $Y$ be a local model of an orbit in the zero fiber and 
$F_Y$ its trivializing homeomorphism defined in Table
\ref{Ta:trivial}. Then on the complement of the exceptional orbits, 
$F_Y$ is a diffeomorphism. 
\end{lemma}

\begin{corollary} \labell{C:alpha-diffeo}
Let $G$ be $\Su(2)$ or $\So(3)$ and $(M, \ow, \Phi)$ be a six
dimensional Hamiltonian $G$-manifold such that $0 \in \Phi(M)$.
Assume that the zero fiber is tall. 
Let $Y$ be a local model of an orbit in the zero fiber and
$F_Y$ its trivializing homeomorphism defined in Table
\ref{Ta:trivial}. Then 
the restriction of $F_Y$ to the reduced space 
on the complement of the exceptional 
orbits is also a diffeomorphism. 
\end{corollary}

\begin{proof}[Proof of Lemma \ref{L:diffeo}]
This is an application of Schwarz \cite{sCh1} and a result of
direct computation. We will prove for the local model $\C^3$ as an
example. 

As in the proof of Lemma \ref{L:homeo}, for $Y=\C^3$, 
the trivializing homeomorphism $F_Y$ is the composition 
$F_2 \circ F_1$, 
where $F_1\colon \C^3 / \So(3) \to \R^3$ is given by
$F_1([q, p]) = (|q|^2, |p|^2, \la q, p \ra)$
and $F_2\colon  \R^3 \to \R_+ \times \C$ is given by 
$F_2(\alpha, \beta, \gamma) =
(\alpha\beta - \gamma^2, \alpha-\beta, 2\gamma)$. 
By Schwarz, $F_1$ pulls back smooth functions to smooth invariant
functions. The inverse $F_2^{-1}$ is given by $F_2^{-1}(a, b, c) = 
(\frac{1}{2}(b+\sqrt{4a+b^2+c^2}), \frac{1}{2}(-b+\sqrt{4a+b^2+c^2},
\frac{c}{2})$. 
By direct computation, $F_2^{-1}$ is smooth except 
when $4a+b^2+c^2 = 0$, i.e., when $a=b=c=0$. So $F_2$ is a
diffeomorphism except when $\alpha = \beta = \gamma =0$. 
So $F_Y$ is a diffeomorphism except at $q=p=0$, the exceptional orbit.

\mute{old proof}{
   Let $E$ denote the set of exceptional orbits. 
   \begin{description}
   \item[Case 1] $Y=\Yone$. $Y \smallsetminus E = Y$. 
	 \[
	\begin{diagram}
   \node{\Yone} \arrow{s} \arrow{se,t}{\tilde{F}}\\
\node{\R^2/S^1 \times \C} \arrow{e,b}{F} \node{\R_+
		   \times \C}
\end{diagram}
\]
   where $F([\mu], z) = (|\mu|^2, z)$. By Schwarz,
   $S^1$-invariant smooth functions are generated by $|\mu|^2$ so
  $F$ is a diffeomorphism. 
	\item[Case 2] $Y=\Ytwo$. Similarly.
 \item[Case 3] $Y=\Ythree$. 
	$F\colon \R^2/S^1 \times (\C \smallsetminus {0})/\Z_2 \to \R_+
	\times \C$ is a diffeomorphism since Schwarz works on $\R^2/S^1$
	follows by the diffeomorphism $z^2$ at $z \neq 0$.
	\item[Case 4] $Y = \C^3$, $F=F_2 \circ F_1$. By Schwarz, $F_1$
	 pulls back smooth functions to smooth invariant functions. $F_2^{-1}$
	 is smooth except $4a+b^2+c^2 = 0$ or $a=b=c=0$. So $F_2$ is a
	diffeomorphism except $\alpha = \beta = \gamma =0$. So $F$ is a
	diffeomorphism except at $p=q=0$, the exceptional orbit. 
	\item[Case 5] $Y = \C^2 \times \C$. $Y \smallsetminus E = Y$. By
	Schwarz, $\Su(2)$-invariant smooth functions are generated by
	$|w|^2$ so $F$ is a diffeomorphism. 
  \end{description}}
\end{proof}

\section{The $\Su(2)$ case}

%Remember what is left to show is that the genus and isotropy data
%determine the quotients up to $\Phi$-diffeomorphism and Condition
%\ref{E:tech} is satisfied.

Let $G$ be $\Su(2)$ and $(M, \ow, \Phi)$ be a six dimensional 
Hamiltonian $G$-manifold with $0 \in \Phi(M)$. 
Assume that the zero fiber $\Phi^{-1}(0)$ is tall. By Table
\ref{Ta:model}, the principal stabilizer of $\Phi^{-1}(0)$ is 
$\Su(2)$ and the 
local model of every orbit in the zero fiber is $Y= \C^2 \times \C$.
By Table \ref{Ta:except}, there is no exceptional orbits. 
Sections \ref{S:qtop} and \ref{S:qsmooth} show that the trivializing
homeomorphism is a diffeomorphism. Near the zero fiber, the quotient
$M/G$ is a smooth manifold with corners since $Y/G$ is
diffeomorphic to $|\Phi(M)|^2 \times \C$. 

Restricting to a smaller neighborhood $V$ of $0 \in \fg^*$, the norm
square of the moment map 
$\overline{\Phi}=|\Phi|^2$ is a proper submersion; there is a
diffeomorphism 
from $\Phi^{-1}(V)/G$ to $(|\Phi(M)|^2 \cap V) \times
(\Phi^{-1}(0)/G)$. 
In particular, $|\Phi(M)|^2 \cap V$ is an interval
and $\Phi^{-1}(V)$ satisfies Condition
\ref{E:tech}.  

The reduced space $\Phi^{-1}(0)/G$ is a Riemann
surface, and is determined by its genus.

Therefore, the genus of the zero fiber determines $\Phi^{-1}(V)$ up
to $\Phi$-diffeomorphisms.  

Shrink $V$ if necessary, by Propositions \ref{P:globallift} and
\ref{P:symplecto}, we have the following version of Theorem A for
$\Su(2)$:  

\begin{proposition}\labell{P:su2case}
Let $G$ be $\Su(2)$. 
Let $(M,\ow, \Phi)$, and $(M', \ow', \Phi')$ be 
compact connected six dimensional Hamiltonian $G$-manifolds such that 
$0 \in \Phi(M) = \Phi'(M')$. Assume that the zero fibers are tall. 
Then there exists an invariant neighborhood $V$ of $0$ in $\fg^*$ over 
which the 
Hamiltonian $G$-manifolds are isomorphic if and only if 
\begin{itemize}
\item their Duistermaat-Heckman functions coincide on $V$, and 
\item their genus at $0$ are the same.
\end{itemize}
\end{proposition}

\section{Grommets}\labell{S:grommet}

Let $G$ be $\So(3)$ and $(M, \ow, \Phi)$ be a six dimensional 
Hamiltonian $G$-manifold with $0 \in \Phi(M)$. Assume that the zero
fiber $\Phi^{-1}(0)$ is tall. According to Table \ref{Ta:model},
there are more than one possible local model for an orbit in the zero
fiber. In particular, some local models have exceptional orbits. 

The arguments in the previous section fail near the exceptional
orbits. Sections 8-11 are dedicated to deal with this problem. The
techniques used here are adapted from \cite{kt1}. 

We start by defining charts at every exceptional orbits so that we
can fix the smooth structure later.

\mute{old def}{
\begin{definition}\labell{D:grommet}
Let $G$ be $\So(3)$ or $\Su(2)$, and let $(M, \ow, \Phi)$ be a
six dimensional Hamiltonian $G$-manifold with $0 \in \Phi(M)$. 
A \textbf{grommet over $\mathbf{0}$} is a $\Phi$-$G$-diffeomorphism 
$\psi\colon  D \to M$ where $D$ is
a $G$-invariant open subset of a local model 
$Y = G \times_H (\fh^0 \times V)$ of an orbit in the zero fiber.   
\end{definition}
}

\begin{definition}\labell{D:grommet}
Let $G$ be a compact Lie group and let $(M, \ow, \Phi)$ be a
complexity one Hamiltonian $G$-manifold with $0 \in \Phi(M)$. 
A \textbf{grommet over $\mathbf{0}$} is a $\Phi$-$G$-diffeomorphism 
$\psi\colon  D \to M$ where $D$ is
a $G$-invariant open subset of a local model 
$Y = G \times_H (\fh^0 \times V)$ of an orbit in the zero fiber.   
\end{definition}

\comment{is this one better?}

\begin{definition}\labell{D:sheet}
Let $G$ be $\So(3)$ or $\Su(2)$, and let $(M, \ow, \Phi)$ be a
six dimensional Hamiltonian $G$-manifold with $0 \in \Phi(M)$. 
Let $Y$ be a local model of an orbit in the zero fiber. 
The \textbf{exceptional sheet} in $Y$ is the subset
\[
S = \{\, [g, \mu, z] \in G \times_H (\fh^0 \times V) \, | \, P_Y([g, \mu, z]) 
= 0\,\}
\]where $P_Y$ is the component function of the trivializing homeomorphism
defined in Table \ref{Ta:trivial}. 
\end{definition}

\begin{remark}
The exceptional orbits are always in the exceptional sheets. However, the
exceptional sheets might include nonexceptional orbits. 
For example, the exceptional sheet of the local model $Y = \C^3$
has a nonexceptional orbit in every $\Phi_Y^{-1}(G\cdot\alpha)$ 
for $\alpha \neq 0 \in \Phi_Y(Y)$.
\end{remark}

\begin{definition}\labell{D:wide}
Let $G$ be $\So(3)$ or $\Su(2)$, and let $(M, \ow, \Phi)$ be a
six dimensional Hamiltonian $G$-manifold with $0 \in \Phi(M)=U$. 
Let $\psi\colon  D \to M$ be a grommet over $0$ where $D$ is
a $G$-invariant open subset of a local model 
$Y = G \times_H (\fh^0 \times V)$ of an orbit in the zero fiber.  
The grommet $\psi$ is \textbf{wide} if $D$ contains its
part of the exceptional sheet. 
That is, $\Phi^{-1}_Y(U) \cap S \subset D$. 
\end{definition}

Wide grommets with pairwise disjoint closures ensure that the
exceptional sheets in different local models are separated.
The following lemma states that 
we can always find such grommets over $0$:

\begin{lemma}\labell{L:widegrommet}
Let $G$ be $\So(3)$ or $\Su(2)$, and let $(N, \ow, \Phi)$ be a
neighborhood of the zero fiber in a six dimensional Hamiltonian
$G$-manifold with $0 \in \Phi(N)$. Assume that the zero fiber is
tall. Let $\{E_j\}$ denote the exceptional orbits in $\Phi^{-1}(0)$.
After replacing $N$ by the preimage of some $G$-invariant
neighborhood of $0$ in $\Phi(N)$, there exist wide grommets
$\psi_j \colon  D_j \to N$ 
over $0$ such that $\psi_j(\{[g, 0, 0]\})= E_j$ and 
$\psi_j(D_j)$ have pairwise disjoint closures.
\end{lemma}

\begin{proof}
By Local Normal Form Theorem \ref{T:normalform}, 
there is a local model $Y_j$ for each
exceptional orbit $E_j \in \Phi^{-1}(0)$ with a moment map 
$\Phi_j \colon  Y_j \to \fg^*$. We can choose
grommets $\psi_j\colon  D_j \subset Y_j \to N$ such that
$\psi_j(\{[g, 0, 0]\}) = E_j$. 

By Lemma \ref{L:homeo} and Definition \ref{D:sheet}, 
the restriction of the norm square of the moment map
$\overline{\Phi}_j$ to $S_j/G$ is a homeomorphism onto its image. 
So there exists a $G$-invariant neighborhood $W_j$ of
$0 \in \fg^*$ such that $S_j \cap D_j = S_j \cap \Phi_j^{-1}(W_j)$.
Let $W = \bigcap W_j$, and replace $N$ by $N \cap \Phi^{-1}(W)$ and $D_j$
by $D_j \cap \Phi^{-1}(W)$. Then the grommets are wide. 

Since each $\psi_j(S_j \cap D_j)$ is closed, 
$\psi_i(S_i \cap D_i) \cap \psi_j(S_j \cap D_j)$ is closed in $N$. It 
does not intersect $\Phi^{-1}(0)$ since the exceptional orbits in the zero
fiber are isolated. Because the moment map is proper, we can find 
a $G$-invariant neighborhood $V \subset W$ of $0$ such that 
$\Phi^{-1}(V)$ does not intersect any of these intersections.

Replace $N$ by $N \cap \Phi^{-1}(V)$ and $D_j$ by $D_j \cap
\Phi^{-1}(V)$. The grommets $\psi$ are still wide and 
$\psi(S_i \cap D_i) \cap \psi_j(S_j \cap D_j)$ is empty. 
Shrinking each $D_j$ to a smaller neighborhood of $S_j \cap D_j$, we
obtain wide grommets with pairwise disjoint closures. 
\end{proof}

\section{Flattening the quotient}\labell{S:flat}

Let $G$ be $\So(3)$ and $(M, \ow, \Phi)$ be a six dimensional 
Hamiltonian $G$-manifold such that $0 \in \Phi(M)$. 
In this section, we show that the quotient $M/G$ near 
the zero fiber is topologically a surface bundle. 
We find a homeomorphism between 
$\Phi^{-1}(V)$ and $(|\Phi(M)|^2 \cap V) \times (\Phi^{-1}(0)/G)$ for
some $G$-invariant neighborhood $V$ of $0$ in $\fg^*$.
Away from the exceptional sheets (see Definition \ref{D:sheet}), 
this homeomorphism is a diffeomorphism. Near the exceptional sheets, 
it is determined by the grommets. 

\comment{away from exceptional sheets or orbits?}

\comment{need to do it also for $\Su(2)$?}

\begin{definition}\labell{D:flatY}
Let $G$ be $\So(3)$ and $(M, \ow, \Phi)$ be a six dimensional
Hamiltonian $G$-manifold such that $0 \in \Phi(M)$. Assume that the
zero fiber is tall. Let $Y$ denote
the local model $G \times_H (\fg^0 \times V)$ of an orbit in the zero
fiber. 
A \textbf{standard flattening} of $Y$ is the map 
\[
\delta\colon  Y /G \to (\text{image}\,\overline{\Phi}_Y) \times
(\Phi^{-1}_Y (0)/G)  
\] given by 
\[
\delta = \left(\overline{\Phi}_Y,\, \overline{P}_0^{-1}
\circ \overline{P}_Y\right)
\]
where $\overline{\Phi}_Y$ is induced by the moment map,
and $\overline{P}_Y\colon  Y/G \to \C$ and
$\overline{P}_0\colon  \Phi^{-1}_Y(0)/G \to \C$ are the component
functions of the trivializing homeomorphisms defined in Section
\ref{S:qtop}, Table \ref{Ta:trivial}.
\end{definition}

By Lemma \ref{L:homeo} and Corollary \ref{C:alpha-homeo}, the
standard flattening is a homeomorphism. By Lemma \ref{L:diffeo} and
Corollary \ref{C:alpha-diffeo}, it is a diffeomorphism away from the
exceptional sheet.

\begin{definition}\labell{D:flat}
Let $G$ be $\So(3)$ and $(N, \ow, \Phi)$ be a neighborhood of the
zero fiber in a six dimensional Hamiltonian $G$-manifold such that 
$0 \in \Phi(N)$. Assume that the
zero fiber is tall. 
The \textbf{flattening} of $N$ about $0$ consists of 
\begin{enumerate}
\item a homeomorphism $\delta \colon  N/G \to
(\text{image}\,\overline{\Phi}) \times (\Phi^{-1}(0)/G)$,  and
\item a wide grommet $\psi_j \colon  D_j \to N$ at each exceptional
orbit $E_j$ in $\Phi^{-1}(0)$ such that $\psi_j(\{[g,0,0]\})=E_j$,
where $D_j$ is 
a $G$-invariant open subset of a local model 
$G \times_H (\fh^0 \times V)$ of an orbit in the zero fiber,
\end{enumerate} 
such that the following two conditions are satisfied:
\begin{enumerate}
\item $\delta$ is a diffeomorphism on the complement of the
exceptional sheets; that is,  
\[
\delta\colon  N/G \smallsetminus \sqcup_j \psi_j(S_j \cap D_j)/G \to
(\text{image}\,\overline{\Phi}) \times (\Phi^{-1}(0)\smallsetminus 
\sqcup_j E_j)/G    
\] is a diffeomorphism.
\item Near the exceptional sheets, $\delta$ is the standard flattening of
the local models; namely, the following diagram commutes:
\[
\begin{diagram}
\node{D_j/G}\arrow{e,t}{\delta_j}\arrow{s,l}{\overline{\psi}_j}
\node{\R_+ \times \left((\Phi^{-1}_j(0) \cap D_j)
/G\right)}\arrow{s,r}{(\text{id}, \overline{\psi}_j)}\\
\node{N/G} \arrow{e,b}{\delta}\node{\R_+ \times
\left({\Phi}^{-1}(0)/G\right)} 
\end{diagram}
\]
where $\overline{\psi}_j\colon D_j/G\to N/G$ is induced by the
grommets, $\Phi_j$ is the moment map on the local model $Y_j \supset
D_j$, and $\delta_j\colon Y_j/G \to (\text{image}\, \overline{\Phi}_j)
\times (\Phi_j^{-1}(0)/G)$ is the standard flattening of $Y_j$. 
\end{enumerate}
\end{definition}

The following proposition asserts that a flattening about $0$ always
exists. 

\begin{proposition}\labell{P:flattening}
Let $G$ be $\So(3)$ and $(M, \ow, \Phi)$ be a six dimensional 
Hamiltonian $G$-manifold such that $0\in \Phi(M)$. Assume that 
the zero fiber is tall. Then there exists a $G$-invariant 
neighborhood $V$ of $0$ in $\Phi(M)$ such that $\Phi^{-1}(V)$ admits
a flattening about $0$. 
\end{proposition} 

\begin{proof}
Let $\{E_j\}$ denote the exceptional orbits in the zero fiber. 
By Lemma~\ref{L:widegrommet}, there exists a $G$-invariant
neighborhood $V$ of $0$ in $\Phi(M)$ and  
wide grommets $\psi_j \colon  D_j \to \Phi^{-1}(V)=N$ 
such that $\psi_j(\{[g, 0, 0]\}) = E_j$ and $\psi_j(D_j)$ have
disjoint closures in $N$. 

The grommets together with the
standard flattening $\delta_j$ of the local models $Y_j$ define a
partial flattening $\delta$ such that the following diagram commutes.
\[
\begin{diagram}
\node{D_j/G}\arrow{s,l}{\overline{\psi}_j}\arrow{e,t}{\delta_j}
\node{\R_+ \times ((\Phi^{-1}_j(0) \cap D_j)/G)}\arrow{s,r}{(\text{id},
\overline{\psi}_j)}\\
\node{\sqcup \overline{\psi}_j(D_j/G)}\arrow{e,b}{\delta}\node{\R_+
%\times ((\Phi^{-1}(0) \cap (\sqcup\psi_j(D_j)))/G)}
\times (\Phi^{-1}(0)/G)}
\end{diagram}
\]

By Lemma \ref{L:diffeo}, $\overline{\Phi}\colon N/G \smallsetminus
\sqcup\psi_j(D_j\cap S_j)/G \to (\text{image}\,\overline{\Phi})$ is a
submersion on the 
complement of the exceptional sheets. The partial flattening $\delta$
defined above then determines an Ehresmann connection on the open set
$\sqcup \psi_j(D_j \smallsetminus S_j)/G$. This connection can be
extended 
on the entire complement of the exceptional sheets with a partition of
unity. Using the parallel transport $\gamma$ associated with the
connection, we define 
$\delta(p) = (\overline{\Phi}(p), \gamma(p))$ for all $p \in N/G$.      
\end{proof}

\comment{this proof should be improved}

This has a few immediate corollaries:

\begin{corollary}\labell{C:sbundle}
Let $G$ be $\So(3)$ and $(M, \ow, \Phi)$ be a six dimensional Hamiltonian
$G$-manifold with $0\in \Phi(M)$. 
Assume that the zero fiber is tall. There exists a $G$-invariant
neighborhood $V$ of $0$ in $\Phi(M)$ 
such that $\Phi^{-1}(V)/G \to \R_+$ is topologically a surface bundle. 
\end{corollary}
 
\begin{corollary}
Let $G$ be $\So(3)$ and $(M, \ow, \Phi)$ be a six dimensional Hamiltonian
$G$-manifold with $0\in \Phi(M)$. Assume that the zero fiber is tall.
Then all the reduced spaces $\Phi^{-1}(G \cdot \alpha)/G$ have the same
genus for $\alpha$ sufficiently close to $0$.
\end{corollary}

\begin{corollary}\labell{C:tech}
Let $G$ be $\So(3)$ and $(M, \ow, \Phi)$ be a six dimensional Hamiltonian
$G$-manifold with $0\in \Phi(M)$. Assume that the zero fiber is tall. 
Then there exists a $G$-invariant neighborhood $V$ of $0$ in $\fg^*$ 
such that the restriction map 
\[
H^*(\Phi^{-1}(V)/G) \to H^*(\Phi^{-1}(G \cdot \alpha)/G)
\]
is an isomorphism for all $\alpha \in V$. In particular,
$\Phi^{-1}(V)$ satisfies Condition \ref{E:tech}. 
\end{corollary}

\section{The associated marked surface}\labell{S:surface}

Let $G$ be $\So(3)$ and $(M, \ow, \Phi)$ be a six dimensional Hamiltonian
$G$-manifold with $0\in \Phi(M)$. Assume that the zero fiber is tall. 
Then the reduced space $\Phi^{-1}(0)$ is topologically a surface. 
We can give it a smooth structure according to the grommets on $M$.
We first define a grommet on a surface. This is simply a notion for a
coordinate chart at a marked point.

\begin{definition}
Let $\Sigma$ denote a closed oriented surface. 
A \textbf{grommet} at a point 
$q \in \Sigma$ is a diffeomorphism $\varphi\colon  B \to
\Sigma$ where $B$ is a neighborhood of $0$ in $\C$ and $\varphi(0) = q$. 
\end{definition}

\begin{definition}
Let $\Sigma$ and $\Sigma'$ be closed oriented surfaces with
labelled marked points
$\{q_i\} \subset \Sigma$ and $\{q_i'\} \subset \Sigma'$. 
Let $\varphi_i$ and $\varphi_i'$ denote the
grommets at these marked points. 
Assume that $q_i$ and $q_i'$ have the same labels for all $i$. 
An orientation preserving diffeomorphism $g\colon  \Sigma \to \Sigma'$ is
\textbf{rigid} if for all $i$
\begin{enumerate}
\item $g(q_i) = q_i'$\,;
\item ${\varphi'_i}^{-1} \circ g \circ \varphi_i$ is
a rotation on some neighborhood of $0 \in \C$.  
\end{enumerate}
\end{definition}

\begin{lemma}\labell{L:rigidmap}
Let $\Sigma$ and $\Sigma'$ be closed oriented surfaces with
labelled marked points
$\{q_i\} \subset \Sigma$ and $\{q_i'\} \subset \Sigma'$. 
Let $\varphi_i$ and $\varphi_i'$ denote the
grommets at these marked points. 
Assume that $q_i$ and $q_i'$ have the same labels for all $i$ and 
that $\Sigma$ and $\Sigma'$ have the same genus. Then 
every bijection from the marked points $\{q_i\}$ to the marked points
with the same labels $\{q_i'\}$ extends to a rigid map
from $\Sigma$ to $\Sigma'$. 
\end{lemma}

\begin{proof}
This uses standard techniques in differential topology. 
For details, see for example \cite{kO}. 
\end{proof}

\begin{remark}\labell{R:chart}
Let $G$ be $\So(3)$ and $(M, \ow, \Phi)$ be a six dimensional Hamiltonian
$G$-manifold. Assume that the zero fiber is tall. Let 
$\psi\colon D \to M$ be a grommet over $0$ where $D$ is a
$G$-invariant open subset of a local model 
$G\times_H(\fh^0 \times V)$ of an orbit in the zero fiber. Assume
that $\psi(\{[g, 0, 0]\}) = \orbit \subset \Phi^{-1}(0)$. 
Then the grommet $\psi$ induces a coordinate chart on the reduced
space $\Phi^{-1}(0)/G$. Explicitly, the map
$\overline{P}_0\colon (D \cap \Phi^{-1}_Y(0))/G \to \C$
given by the trivializing homeomorphism in Table \ref{Ta:trivial} is a
homeomorphism onto its image $B$. And 
$\varphi=\overline{\psi}\circ \overline{P}_0^{-1}\colon B \to
\Phi^{-1}(0)/G$ is a homeomorphism onto its image such that
$\varphi(0) = \orbit/G$, where $\overline{\psi}\colon D/G \to M/G$ is
induced by $\psi$.   
\end{remark}

\comment{do i need to indicate somewhere before that
$\overline{P}_Y([g, 0, 0])=0$?}

\begin{definition} \labell{D:surface}
Let $G$ be $\So(3)$ and $(N, \ow, \Phi)$ be a $G$-invariant
neighborhood of the zero fiber in a six dimensional Hamiltonian
$G$-manifold. Assume that the zero fiber is tall. For each
exceptional orbit $E_j$ in $\Phi^{-1}(0)$, let $\psi_j\colon D_j \to
N$ be the grommet over $0$ such that $\psi_j(\{[g, 0, 0]\}) = E_j$.  
The \textbf{associated marked surface} of $N$ consists the following
data: 
\begin{enumerate}
\item The connected oriented topological surface $\Sigma = \Phi^{-1}(0)/G$.
\item The set of marked points $\{q_j\} \in \Sigma$ corresponding 
to the exceptional orbits $\{E_j\} \in \Phi^{-1}(0)$; i.e. for each $j$,
$q_j=E_j/G$.  
\item The smooth structure on $\Sigma$ given by the following
coordinate charts. For each exceptional orbit $E_j$ in
$\Phi^{-1}(0)$, take the given grommet. For each nonexceptional orbit
$\orbit$ in $\Phi^{-1}(0)$, choose an arbitrary grommet such that
$\psi(\{[g, 0, 0]\})=\orbit$. For each grommet, take the induced
coordinate chart as described in Remark \ref{R:chart}. 
\item The grommets on $\Sigma$ at the marked points $\{q_j\}$ given
by the above coordinate charts. 
\item A label at each marked point $q_j$ describing the isotropy 
representation of the corresponding exceptional orbit $E_j$.  
\end{enumerate} 
\end{definition}

\comment{need to specify the principal stabilizer later on in the
statements since it is not specified in the associated surface}

\comment{Is it better to define associated surface 
for all $\alpha$ or is it necessary? It will need a definition of
grommets over $\alpha$}

\section{Diffeomorphism between quotients}

Here we show that a diffeomorphism between the associated marked surfaces 
that behaves nicely near the marked points 
induces a $\Phi$-diffeomorphism between the quotients. 

\begin{proposition} \labell{P:rigidext}
Let $G$ be $\So(3)$ and let $(N, \ow, \Phi)$ and $(N', \ow', \Phi')$
be neighborhoods of the zero fibers in six dimensional Hamiltonian
$G$-manifolds such that $0 \in \Phi(N) = \Phi'(N')$. 
Assume that the zero fibers are tall and have the same principal
stabilizer. Assume that $N$ and $N'$ admit flattenings about
$0$. 
Let $\Sigma$ and $\Sigma'$ denote the associated marked surfaces of
$N$ and $N'$ respectively. 
Then any rigid map $h \colon  \Sigma \to \Sigma'$ extends to a
$\Phi$-diffeomorphism $H\colon  N/G \to N'/G$.
\end{proposition}

\begin{proof}
If there exists a rigid map $h\colon \Sigma \to \Sigma'$, the labels
of the marked points on $\Sigma$ and $\Sigma'$ are the same.
Topologically, $\Phi^{-1}(0)/G=\Sigma$ and
${\Phi'}^{-1}(0)/G=\Sigma'$, so the isotropy data at $0$ are the same
for $N$ and $N'$. Let $\overline{\Phi}$ and $\overline{\Phi'}$ be 
the norm squares of the moment maps. We know that 
$\text{image}\,\overline{\Phi} = \text{image}\,\overline{\Phi'}$.

Let $\delta$ and $\delta'$ denote the homeomorphisms given by the
flattenings of $N$ and $N'$ about $0$.  
We define $H = {\delta'}^{-1} \circ (\text{id}, h) \circ \delta$ so that
the following diagram commutes:
\[
\begin{diagram}
\node{N/G}\arrow{e,t}{H}\arrow{s,l}{\delta}\node{N'/G}\arrow{s,r}{\delta'}\\
\node{(\text{image}\,\overline{\Phi}) \times \Sigma}\arrow{e,b}{(\text{id},
h)}\node{(\text{image}\,\overline{\Phi'}) \times \Sigma'}
\end{diagram}
\]
We claim that $H$ is a $\Phi$-diffeomorphism. 

It is orientation
preserving since $h$ is. $H^*\overline{\Phi'}= \overline{\Phi}$ since 
$H$ is induced from an identity map between
$(\text{image}\,\overline{\Phi})$ and 
$(\text{image}\,\overline{\Phi'})$. We need to prove that it is a
diffeomorphism and it locally lifts in a neighborhood of 
each exceptional orbit.

Let $\{E_j\}$ and $\{E_j'\}$ denote the exceptional orbits in
$\Phi^{-1}(0)$ and ${\Phi'}^{-1}(0)$, respectively.  
The rigid map $h$ determines an identification between the
exceptional orbits $E_j$ and $E_j'$ 
with the same isotropy representation.
We can then denote the local model by $Y_j$ for both 
$E_j \in \Phi^{-1}(0)$ and $E_j'\in {\Phi'}^{-1}(0)$. 
Let $$\psi_j\colon D_j \to N \qquad \text{and} \qquad \psi_j'\colon
D_j'\to N'$$ 
denote the grommets given by the flattenings of $N$ and $N'$ about
$0$, where $D_j \subset Y_j$ and $D_j' \subset Y_j$ are $G$-invariant
open subsets. Let $S_j$ and $S_j'$ denote the exceptional sheets. 
     
Since the homeomorphisms $\delta$ and $\delta'$ given in the
flattenings are diffeomorphisms on the complement of the exceptional
sheets, and since the smooth structures on the reduced space 
$\Phi^{-1}(0)/G$ and the associated surface $\Sigma$ agree
off the exceptional orbits, the restriction of $H$ to
$({N\smallsetminus \sqcup_j \psi_j(S_j \cap D_j)})/{G}$ is a
diffeomorphism. This is easy to see from the following diagram, where
$$\Tilde{h}\colon \left({\Phi^{-1}(0)\smallsetminus  \sqcup_j E_j}\right)/{G}
\to \left({{\Phi'}^{-1}(0)\smallsetminus  \sqcup_j E'_j}\right)/{G}$$ is a
diffeomorphism induced by $h\colon \Sigma \to \Sigma'$: 
\[
\begin{diagram}
\node{\frac{N\smallsetminus \sqcup_j \psi_j(S_j \cap D_j)}{G}}\arrow{e,t}{H}
\arrow{s,l}{\delta}\node{\frac{N'\smallsetminus \sqcup_j \psi'_j(S'_j \cap
D'_j)}{G}}\arrow{s,r}{\delta'}\\
\node{(\text{image}\,\overline{\Phi})\times
\frac{\Phi^{-1}(0)\smallsetminus \sqcup_j 
E_j}{G}}\arrow{r,b}{(id, \Tilde{h})}\node{(\text{image}\,\overline{\Phi'})\times
\frac{{\Phi'}^{-1}(0)\smallsetminus  \sqcup_j E'_j}{G}}
\end{diagram}
\] 

It remains to show that $H$ is a $\Phi$-diffeomorphism in a
neighborhood of each exceptional sheets $\psi_j(S_j \cap D_j)/G$. 

Let $\varphi_j\colon B_j \to \Sigma$ and $\varphi'_j\colon B_j' \to
\Sigma'$ be the grommets of the associated marked surfaces.  
The fact that $h$ is
rigid implies that ${\varphi'_j}^{-1} \circ h \circ \varphi_j$ is a rotation by
$a_j \in S^1$ on some neighborhood of $0$. That is, we have the
following diagram:
\[
\begin{diagram}
\node{B_j}\arrow{e,t}{w \mapsto a_j
w}\arrow{s,l}{\varphi_j}\node{B_j'}\arrow{s,r}{\varphi'_j}\\
\node{\Sigma}\arrow{e,b}{h}\node{\Sigma'} 
\end{diagram}
\]
For each local model $Y_j$ of an exceptional orbit, we define a
map $f_j\colon Y_j \to Y_j$ as follows: 
	\begin{description}
	\item[Case 1] $Y_j = \So(3) \times_{N(S^1)}(\R^2 \times \C), \quad
	f_j([g, \mu, z]) = [g, \mu, \pm a_j^{1/2}z]$.
	\item[Case 2] $Y_j = \C^3, \quad f_j(x, y, z) = \pm (a_j^{1/2}x,\,
	a_j^{1/2}y,\, a_j^{1/2}z)$.
	\end{description}
This map $f_j$ is an equivariant symplectomorphism which respects the
moment maps. It induces a $\Phi$-diffeomorphism $\Tilde{f}_j\colon
Y_j/G \to Y_j/G$ on the quotient. 
In some neighborhood of $\psi_j(S_j\cap D_j)$, if we identify
$Y_j/G$ with $(\text{image}\,{\overline{\Phi}_j}) \times \C$ by the
trivializing homeomorphism, $\Tilde{f}_j$ sends $(\alpha, w)$ to
$(\alpha, a_j w)$, and therefore agrees with $H$.   
\end{proof}

\comment{where is the principal stabilizer of the zero fiber?}

\section{Proof of the local uniqueness theorem over zero}\labell{S:part1end}

We now arrive at the stage to prove our main local theorem.

\begin{theorema}[Local Uniqueness over $0$]%\labell{T:locun}
Let $G$ be $\Su(2)$ or $\So(3)$. 
Let $(M,\ow, \Phi)$ and $(M', \ow', \Phi')$ be 
compact connected six dimensional Hamiltonian $G$-manifolds such that 
$0 \in \Phi(M) = \Phi'(M')$. 
There exists a $G$-invariant neighborhood $V$ of $0$ in $\fg^*$ over which the 
Hamiltonian $G$-manifolds are isomorphic if and only if 
\begin{itemize}
\item their Duistermaat-Heckman functions coincide on $V$, 
%(in $\So(3)$ case, they only have to agree at some $\alpha \in V$), 
\item their isotropy data and genus at $0$ are the same,
\item their principal stabilizers of the zero fibers are the same,
\item if the zero fibers are tall with principal stabilizer $S^1$, 
the first Stiefel-Whitney classes of $\Phi^{-1}(0)$ and
${\Phi'}^{-1}(0)$ in $H^1(M_0^{\text{reg}}; \Z_2)$ and
$H^1({M'}_0^{\text{reg}}; \Z_2)$
are equal (under an identification of $M_0$ and $M_0'$ that respects
the isotropy data). 
\end{itemize}
\end{theorema}

\begin{proof}
It is clear that these conditions are necessary conditions when
applicable. We show that they are also sufficient conditions.
 
If the zero fibers are short, this is Proposition
\ref{P:shortzerofiber}. 

If $G$ is $\Su(2)$ and the zero fibers are tall, this is Proposition
\ref{P:su2case}. 

Assume that $G$ is $\So(3)$ and that the zero fibers are tall.  
By Proposition \ref{P:flattening}, there exists a $G$-invariant
neighborhood $V$ of $0$ in $\Phi(M)=\Phi'(M')$ such that 
$\Phi^{-1}(V)$ and ${\Phi'}^{-1}(V)$ admit flattenings about $0$. 
By assumption, their genus and isotropy data are the same at $0$.  
Then the associated marked surfaces $\Sigma$ and $\Sigma'$ have the
same genus and labels. By Lemma \ref{L:rigidmap}, there exists a rigid
map $h\colon  \Sigma \to \Sigma'$. By Proposition \ref{P:rigidext}, it
extends to a $\Phi$-diffeomorphism $\psi\colon \Phi^{-1}(V)/G \to
{\Phi'}^{-1}(V)/G$. Now Condition \ref{E:tech} is satisfied because
of the flattenings. By assumption, the Duistermaat-Heckman functions
of $\Phi^{-1}(V)$ and ${\Phi'}^{-1}(V)$ are the same and the zero
fibers have the same principal stabilizer and the first
Stiefel-Whitney class. 
Proposition \ref{P:globallift} implies that there exists a
$\Phi$-$G$-diffeomorphism from $\Phi^{-1}(V)$ to ${\Phi'}^{-1}(V)$,
and Proposition \ref{P:symplecto} implies that there exists an
equivariant symplectomorphism from $\Phi^{-1}(V)$ to
${\Phi'}^{-1}(V)$ which respects the moment maps.  
\end{proof}

\section{Symplectic cross-section}\labell{S:part2begin}

In this section we begin to study the preimage of a neighborhood away
from zero. 
First, we factor out the coadjoint orbit directions in the
sense of the symplectic cross-section theorem (cf. Theorem 26.7 in
\cite{gs2}).   

\begin{theorem}
Let $G$ be a compact Lie group and let $(M, \ow, \Phi)$ be 
a Hamiltonian $G$-manifold. Suppose that 
$S$ is a submanifold of
$\fg^*$ passing through a point $\alpha \in \fg^*$ satisfying 
$T_\alpha S \oplus T_\alpha (G\cdot\alpha) = \fg^*$ 
and suppose that $S$ is $G_\alpha$-invariant.
Then for a small enough $G_\alpha$-invariant neighborhood $B$ of
$\alpha$ in $S$ the preimage $\Phi^{-1}(B)$ is a symplectic
submanifold of $(M, \ow)$ and the action of $G_\alpha$ on $\Phi^{-1}(B)$ is
Hamiltonian. Its moment map is the restriction of $\Phi$ on
$\Phi^{-1}(B)$ followed
by the projection onto $T_\alpha S \simeq \fg_\alpha^*$, the dual of
the Lie algebra of the stabilizer of $\alpha$.
\end{theorem}

The submanifold $X = \Phi^{-1}(B)$ is called a \textbf{symplectic
cross-section}. 
It is proved in \cite{gls} that if $G$ is compact we can choose 
the manifold $S$ and $B$ be so large that it is all the interior 
of the Weyl chamber. The set $\Phi^{-1}(G \cdot B) = G \cdot
\Phi^{-1}(B) = G \cdot X$ is an open subset of $M$,
which is $G$-equivariantly isomorphic to the associated bundle $G
\times_{G_\alpha}X$, and the map $\pi\colon  G \times_{G_\alpha} X \to G \cdot
\alpha$ given by $[g, x] \mapsto g \cdot \alpha$ is a symplectic fibration. The
symplectic connection on the bundle $G \times_{G_\alpha} X \to G
\cdot \alpha$ is the same as the connection 
determined by the splitting $T_\alpha S \oplus T_\alpha
(G\cdot\alpha) = \fg^*$. While the symplectic form on $G
\times_{G_\alpha}X$ comes naturally from its identification with 
$G \cdot X \hookrightarrow M$, it can be reconstructed from its
restriction to the fiber $X$ and the connection corresponding to the
splitting.   

Let $G$ be $\Su(2)$ or $\So(3)$, and let $(M, \ow, \Phi)$ be a six
dimensional Hamiltonian $G$-manifold. For all $\alpha \neq 0 \in
\fg^*$, we have $G_\alpha = S^1$ and 
$G\cdot \alpha = S^2$. 
The symplectic cross-section $X$ is the
preimage of 
an open ray $\R_{>0} = \R_+ \smallsetminus \{0\}$. 
It is symplectic, connected, four dimensional and it has 
a Hamiltonian circle action. Its moment map 
$\Phi_X \colon  X \to \R_{>0}$ is given by
\[
g \cdot \Phi_X(x) = \Phi(\,[\,g, x\,]\,)
\] where $\Phi$ is the moment map 
of the $G$ action on $G \times_{G_\alpha} X$. 
Let $A$ denote the connection one-form. 
The symplectic form on $G \times_{G_\alpha} X$ can 
be reconstructed by
\begin{equation}\labell{E:coupling}
\ow = \ow_X - d \la \Phi_X, A \ra + \pi^* \ow_{S^2}
\end{equation}
where $\ow_X$ is a symplectic form on the symplectic cross-section
$X$, and $\pi^*\ow_{S^2}$ is the pull-back of the natural
symplectic form on the coadjoint orbit $G \cdot \alpha = S^2$. 

\comment{the notation $\Phi_X$ is not good} 

A special case of Theorem B follows from the symplectic cross-section
theorem:

\begin{proposition}\labell{P:nonzero}
Let $G$ be $\Su(2)$ or $\So(3)$. 
Let $(M, \ow, \Phi)$ and $(M', \ow',\Phi')$ be compact connected
six dimensional Hamiltonian $G$-manifolds such that 
$0 \notin \Phi(M)=\Phi'(M)$. Then $M$ and $M'$ are isomorphic if and
only if they have the same principal stabilizer, the same
Duistermaat-Heckman function, the same genus, and the same 
isotropy skeleton. 
\end{proposition}

\begin{proof}
Since $M$ and $M'$ have the same principal stabilizer
$\Z_n \simeq \Z/n\Z$,  their symplectic cross-sections $X$ and $X'$ inherit 
actions of $S^1/\Z_n \simeq S^1$ from the actions of 
$G$ on $M$ and $M'$. The symplectic
cross-sections are then four dimensional Hamiltonian $S^1$-manifolds
on which $\Z_n$ acts trivially.   
By \cite{k1} or \cite{kt2}, these are determined by the 
Duistermaat-Heckman function, the
genus, and the isotropy skeleton. Hence $X$ and $X'$ are
isomorphic and $M \simeq G \times_{S^1} X$ and $M' \simeq G
\times_{S^1} X'$ are isomorphic. 
\end{proof}

We now examine some previous definitions and properties on a
Hamiltonian $G$-manifold in terms of its symplectic 
cross-section.   

\begin{lemma} \labell{L:phig121}
Let $G$ denote $\Su(2)$ or $\So(3)$.
Let $(M, \ow, \Phi)$ and $(M', \ow', \Phi')$ be  
Hamiltonian $G$-manifolds and let 
$X$ and $X'$ be their symplectic cross-sections 
$\Phi^{-1}(\R_{>0})$ and ${\Phi'}^{-1}(\R_{>0})$ respectively.
Identify $M \smallsetminus \Phi^{-1}(0)$ and $M' \smallsetminus
{\Phi'}^{-1}(0)$ with $G \times_{S^1} X$ and $G \times_{S^1} X'$.  
There exists a one-to-one correspondence between 
$\Phi$-$G$-diffeomorphisms 
$F \colon  G \times_{S^1} X \to G \times_{S^1} X'$ and
$\Phi_X$-$S^1$-diffeomorphisms $f\colon  X \to X'$.
\end{lemma}

\begin{proof}
One direction is easy. Assume that there exists a
$\Phi_X$-$S^1$-diffeomorphism $f\colon  X \to X'$. The map 
$F \colon  G \times_{S^1} X \to G \times_{S^1} X'$ defined by 
$F([g,x]) = [g, f(x)]$ is a $\Phi$-$G$-diffeomorphism. 

Now assume there exists a $\Phi$-$G$-diffeomorphism $F\colon  G
\times_{H} X \to G \times_{H} X'$. The fact that $F$ is
equivariant implies
that $F([g,x]) = [gf_1(x), f_2(x)]$ for some smooth
functions $f_1\colon  X \to G$ and $f_2\colon  X \to X'$.
Since $F$ respects the moment maps, $g\cdot\Phi_X(x) =
gf_1(x)\cdot\Phi_X'(f_2(x))$. Note $\Phi_X$ and $\Phi_X'$ map 
into $\R_{>0}$. We derive
that $f_1(x) \in S^1$ and $\Phi_X(x) = \Phi_X'(f_2(x))$. We rewrite 
$F([g, x]) = [gf_1(x), f_2(x)] = [g, f_1(x)f_2(x)]= [g, f(x)]$. 

To show that $f$ is a $\Phi_X$-$S^1$-diffeomorphism, we need to show
that $f$ is an orientation preserving diffeomorphism which respects
the moment maps.  
Since $F$ is well-defined under the
$S^1$ action on $G \times X$, we have $F([g, x]) = [g, f(x)] = 
F([ga^{-1}, ax]) = [ga^{-1},f(ax)]$ for all $a \in S^1$. 
So $f(ax) = a f(x)$ for all $a \in S^1$; namely, 
$f$ is $S^1$-equivariant. We know that this
function $f\colon X \to X'$ also respects the
moment map since $\Phi_X'$ is $S^1$-invariant and 
$\Phi_X'(f(x)) = \Phi_X'(f_1(x)f_2(x)) = \Phi_X'(f_2(x)) = \Phi_X(x)$. 
Finally, the fact that 
$F$ preserves the orientation on $G \times_{S^1} X$ implies
that $f$ preserves the orientation on $X$.  
\end{proof}

\mute{omit}{
Let $G$ be a group, $H$ a subgroup of $G$, and $X$ an $H$-space. The
stabilizer $G_{[g, x]}$ at $[g, x]$ of the induced $G$-space $G
\times_H X$ is given by $G_{[g, x]} = gH_x g^{-1}$. Therefore, we
have a correspondence between the exceptional orbits of a Hamiltonian
$G$-manifold and those of its symplectic cross-section:

\begin{lemma}
Let $G$ be $\Su(2)$ or $\So(3)$ and $(M, \ow, \Phi)$ be a 
Hamiltonian $G$-manifold. 
Let $X$ denote the symplectic cross-section 
$\Phi^{-1}(\R_{>0})$. 
The $S^1$-orbit $\orbit$ is exceptional in
$\Phi^{-1}_X(\alpha)$ if 
and only if the $G$-orbit $G\cdot\orbit$ is exceptional in
$\Phi^{-1}(G\cdot\alpha)$.    
\end{lemma}
}

\begin{definition}
Let $\R$ denote the dual of the Lie algebra of $S^1$. 
Let $M$ and $M'$ be oriented manifolds with $S^1$ actions and
$S^1$-equivariant maps $\Phi\colon M \to \R$ and 
$\Phi'\colon  M' \to \R$. 
A \textbf{$\Phi$-diffeomorphism} from $M/S^1$ to $M'/S^1$ 
is an orientation preserving diffeomorphism $\psi\colon  M/S^1 \to
M'/S^1$ such that $\psi^* \overline{\Phi'} = \overline{\Phi}$ and
such that $\psi$ and $\psi^{-1}$ lift to $\Phi$-$S^1$-diffeomorphisms
in a neighborhood of each exceptional orbit. Here $\overline{\Phi}$
and $\overline{\Phi'}$ are induced by the moment maps as in
(\ref{E:norm2}). 
\end{definition}

\begin{lemma}\labell{L:phi121}
Let $G$ denote $\Su(2)$ or $\So(3)$. 
Let $(M, \ow, \Phi)$ and $(M', \ow', \Phi')$ be 
Hamiltonian $G$-manifolds and 
let $X$ and $X'$ be their symplectic cross-sections 
$\Phi^{-1}(\R_{>0})$ and ${\Phi'}^{-1}(\R_{>0})$ respectively. 
Identify $M \smallsetminus \Phi^{-1}(0)$ and 
$M' \smallsetminus {\Phi'}^{-1}(0)$ with
$G\times_{S^1} X$ and $G \times_{S^1} X'$.
There exists a one-to-one correspondence between 
$\Phi$-diffeomorphisms $F\colon  (G \times_{S^1} X)/G \to (G
\times_{S^1} X')/G$ and 
$\Phi_X$-diffeomorphisms $f\colon  X/S^1 \to X'/S^1$.
\end{lemma}

\begin{proof}
Let $H$ denote $S^1$.
Consider the map $\overline{i}\colon X/H \to (G \times_H X)/G$ induced by
the inclusion $i\colon X \to G \times_H X$ such that $i(x) = [e, x]$. It is
easy to check that $\overline{i}$ is a homeomorphism. We have the
commutative diagram that gives the correspondence between a
$\Phi$-diffeomorphism $F$ and a $\Phi_X$-diffeomorphism $f$. 
\[
\begin{diagram}
\node{X/H}\arrow{e,t}{f}\arrow{s,r}{\overline{i}}
\node{X'/H}\arrow{e,t}{\overline{\Phi}_X'}\arrow{s,r}{\overline{i}}
\node{\R_{>0}}\arrow{s,r}{\text{id}}\\
\node{(G \times_H X)/G}\arrow{e,b}{F}
\node{(G \times_H X')/G}\arrow{e,b}{\overline{\Phi}'}
\node{\R_{>0}}
\end{diagram}
\]

Let $H_x$ denote the stabilizer of a point $x \in X$ with respect to
the $H$ action on $X$. The stabilizer of $[g, x] \in G \times_H
X$ with respect to the $G$ action on $G \times_H X$ is given by 
$G_{[g, x]} = gH_xg^{-1}$.  
Since $\Phi([g,x]) = g \cdot \Phi_X (x)$, 
there is a one-to-one correspondence between exceptional $S^1$-orbits 
in $\Phi_X^{-1}(\alpha)$ and exceptional $G$-orbits in
$\Phi^{-1}( G \cdot \alpha)$. By Lemma \ref{L:phig121}, 
if the $\Phi_X$-diffeomorphism $f$ lifts to a
$\Phi_X$-$S^1$-diffeomorphism 
in a neighborhood of an exceptional $S^1$-orbit, the induced
$\Phi$-diffeomorphism $F = \overline{i} \circ f \circ
\overline{i}^{-1}$ lifts to a $\Phi$-$G$-diffeomorphism in a
neighborhood of the corresponding exceptional $G$-orbit.
\end{proof}

\begin{lemma}\labell{L:lift121}
Following the notations in Lemma \ref{L:phig121} and
\ref{L:phi121}, a $\Phi$-diffeomorphism $F$ lifts to a
$\Phi$-$G$-diffeomorphism $\Tilde{F}$ if and only if the
corresponding $\Phi_X$-diffeomorphism $f$ lifts to the corresponding
$\Phi_X$-$S^1$-diffeomorphism $\Tilde{f}$.
\end{lemma}

\begin{proof}
Let $[\ ,\ ]$ denote the $S^1$ equivalence class and let $\la\ , \
\ra$ denote the $G$ equivalence class. 
A $\Phi$-diffeomorphism $F \colon  (G\times_{S^1}X)/G \to
(G\times_{S^1}X')/G$ lifts to a $\Phi$-$G$-diffeomorphism
$\Tilde{F} \colon  G\times_{S^1}X \to G\times_{S^1}X'$ 
if and only if
\begin{equation}\labell{E:1}
\la\Tilde{F}([g,x])\ra = F(\la[g,x]\ra),
\end{equation}
and a $\Phi_{X}$-diffeomorphism $f \colon  X/S^1 \to X'/S^1$ lifts to a
$\Phi_{X}$-$S^1$-diffeomorphism $\Tilde{f}\colon X \to X'$ if and only
if 
\begin{equation}\labell{E:2}
[\Tilde{f}(x)]=f([x]).
\end{equation}
By Lemma \ref{L:phig121}, 
\begin{equation}\labell{E:3}
\Tilde{F}([g,x]) = [g, \Tilde{f}(x)].
\end{equation}
By Lemma \ref{L:phi121}, 
\begin{equation}\labell{E:4}
F(\la[e,x]\ra) = \la[e, f([x])]\ra.
\end{equation}
We only need to show that (\ref{E:2}), (\ref{E:3}), and (\ref{E:4}) together
imply (\ref{E:1}), but this is easy.
\[
F(\la[g, x]\ra) = F(\la[e,x]\ra) = \la[e, f([x])]\ra =
\la[e, [\Tilde{f}(x)]]\ra
\]
clearly equals
\[
\begin{split}
\la\Tilde{F}([g, x])\ra &= \la[g, \Tilde{f}(x)]\ra = 
\la[e,\Tilde{f}(x)]\ra = \la[gh^{-1}, h\Tilde{f}(x)]\ra \\
&= \la[e,h\Tilde{f}(x)]\ra = \la[e, [\Tilde{f}(x)]]\ra
\end{split}
\] where $h$ is some element in $S^1$.
\end{proof}

The moment image of a Hamiltonian $S^1$-manifold can be
translated by a constant $\alpha \in {\mathfrak{s}}^*\simeq \R$.
Hence, the Local Normal Form Theorem \ref{T:normalform} applies 
for any orbit in any $\alpha$ fiber. 
It follows from \cite{kt1} 
that trivializing homeomorphisms exist, and 
that grommets, exceptional sheets, wide grommets, 
flattenings, and associated marked surfaces described in Sections
\ref{S:grommet}, \ref{S:flat}, and \ref{S:surface} 
are defined over any tall $\alpha$ fiber for the symplectic 
cross-section $X$ of a six dimensional Hamiltonian $G$-manifold
where $G=\Su(2)$ or $\So(3)$.

Let $Y$ be a local model for an $S^1$-orbit $\orbit$ in 
$\Phi^{-1}_X(\alpha)$ for $\alpha \in \R_{>0}$. 
The associated bundle $G\times_{S^1}Y$ is a local model for the $G$-orbit
$G\cdot \orbit$ in $\Phi^{-1}(G\cdot\alpha)$. 
A trivializing homeomorphism on $Y$
determines a trivializing homeomorphism on $G\times_{S^1} Y$. 
We can then define grommets, exceptional sheets, wide grommets,
flattenings, and associated marked surfaces of $M$ away from $0$
in a similar fashion. 

\comment{need to elaborate}

The following propositions and lemmas are derived from \cite{kt1}
based on the properties of the symplectic cross-section.

%\begin{definition} 
%flattening
%\end{definition}

\begin{proposition}\labell{P:alphaflat}
Let $G$ be $\Su(2)$ or $\So(3)$ and let
$(M, \ow, \Phi)$ be a six dimensional Hamiltonian
$G$-manifold. Assume that $\Phi^{-1}(\alpha)$ is tall for $\alpha \in
\Phi(M) \cap \R_{>0}$. Then there
exists a neighborhood $I$ of $\alpha$ in $\Phi(M) \cap
\R_{>0}$ such that the preimage $\Phi^{-1}(G \cdot I)$ 
admits a flattening about $\alpha$. 
\end{proposition}

%This is Proposition 9.5 in \cite{kt1}.

%\begin{definition}
%associated marked surface
%\end{definition}

\begin{proposition}\labell{P:alpharigidext}
Let $G$ be $\Su(2)$ or $\So(3)$. 
Let $(M, \ow, \Phi)$ and $(M', \ow', \Phi')$ be six dimensional 
Hamiltonian $G$-manifolds equipped with flattenings about $\alpha
\in \Phi(M)=\Phi'(M')$. Let $\Sigma$ and $\Sigma'$ be the associated
marked surfaces of $M$ and $M'$, respectively. 
Then any rigid map $h\colon \Sigma \to \Sigma'$
extends to a $\Phi$-diffeomorphism $g\colon M/G \to M'/G$. 
\end{proposition}

%This is Proposition 11.1 in \cite{kt1}.

\begin{proposition}
Let $G$ be $\Su(2)$ or $\So(3)$. Let $(M, \ow, \Phi)$ and $(M', \ow',
\Phi')$ be compact connected six dimensional Hamiltonian
$G$-manifolds. 
Assume that their $\alpha$ fibers are tall for $\alpha \neq 0$ in
$\Phi(M) = {\Phi'}(M')$.  Then there exists a $G$-invariant 
neighborhood $V$ of $\alpha$ such that $\Phi^{-1}(V)$ and 
${\Phi'}^{-1}(V)$ are $\Phi$-diffeomorphic if and only if the
reduced spaces $\Phi^{-1}(G\cdot\alpha)/G$ and
${\Phi'}^{-1}(G\cdot\alpha)/G$ have the same 
isotropy data and genus. 
\end{proposition}

\begin{proof}
%Let $X$ and $X'$ denote the symplectic cross-sections
%$\Phi^{-1}(\R_{>0})$ and ${\Phi'}^{-1}(\R_{>0})$. 
%We know that $X$ and $X'$ are four dimensional
%Hamiltonian $S^1$-manifolds. Let $\Phi_X$ and $\Phi'_X$ denote the
%moment maps on $X$ and $X'$, respectively.   
Without loss of generality, we can assume $\alpha \in \R_{>0}$. 
By Proposition \ref{P:alphaflat}, there exists a neighborhood $I$ of
$\alpha$ in $\R_{>0}$ such that $\Phi^{-1}(G\cdot I)$ and
${\Phi'}^{-1}(G\cdot I)$ admit flattenings about $\alpha$. Since
the isotropy data and genus of the
reduced spaces $M_\alpha$ and $M'_\alpha$ are the same, 
the associated marked surfaces $\Sigma$ and $\Sigma'$ 
have the same genus and labels. By
Lemma \ref{L:rigidmap}, there exists a rigid map $h\colon \Sigma \to
\Sigma'$. By Proposition \ref{P:alpharigidext}, there exists a
$\Phi$-diffeomorphism $g\colon \Phi^{-1}(G\cdot I)/G \to
{\Phi'}^{-1}(G\cdot I)/G$. 
\end{proof}

\comment{need to explain the identification between $X$ and $X
\subset M$ and also $\Phi_X$ and $\Phi|_X$}

\comment{such a mess!}

\section{Global structure of the orbit space $M/G$}

Let $G$ be $\Su(2)$ or $\So(3)$ and let $(M, \ow, \Phi)$ be a six
dimensional Hamiltonian $G$-manifold. 
Assume that every moment fiber is tall. Corollary \ref{C:sbundle} 
states that there exists a $G$-invariant 
neighborhood $V$ of zero such that $\Phi^{-1}(V)/G$ is 
topologically a trivial surface bundle. Most importantly, we have
shown that the restriction map 
$H^i(\Phi^{-1}(V)/G; \Z) \to H^i( \Phi^{-1}(G \cdot\alpha)/G; \Z)$
is one-to-one for $i=1, 2$ and $\alpha \in V$. 

Assume $\alpha \neq 0$ and $\alpha \in \Phi(M)$. 
There also exists a
$G$-invariant neighborhood $V$ of $\alpha$ such that 
the restriction map 
$H^i(\Phi^{-1}(V)/G; \Z) \to H^i( \Phi^{-1}(G \cdot\beta)/G; \Z)$
is one-to-one for $i=1, 2$, and $\beta \in V$. 
This essentially follows from \cite{kt1} 
since $\Phi^{-1}(G\cdot\beta)/G = \Phi_X^{-1}(\beta)/S^1$, and 
$\Phi^{-1}(G\cdot I)/G=(G \times_{S^1}
\Phi_X^{-1}(I))/G=\Phi_X^{-1}(I)/S^1$ for any interval $I \subset
\R_{>0}$.
Here $X$ denotes the symplectic cross-section
$\Phi^{-1}(\R_{>0})$ and $\Phi_X$ is the corresponding moment map for
the Hamiltonian circle action on $X$.

\comment{this is flattening, quote specific from KT1?}

In this section, we will show that the injectivity of the restriction
map holds not only locally but also for the manifold $M$.

\begin{proposition}\labell{P:techM}
Let $G = \Su(2)$ or $\So(3)$ and let $(M, \ow, \Phi)$ be a connected
six dimensional Hamiltonian $G$-manifold. Then the restriction map 
$H^i(M/G; \Z) \to H^i( \Phi^{-1}(G \cdot\alpha)/G; \Z)$ 
is one-to-one for any $\alpha \in \Phi(M)$ and $i = 1, 2$. 
\end{proposition}

We need several observations to carry out the proof.

\begin{lemma}\labell{L:genus0}
Let $G$ be $\Su(2)$ or $\So(3)$ and 
let $(M, \ow, \Phi)$ be a six dimensional Hamiltonian $G$-manifold
with a moment map $\Phi\colon M \to \R^3$. Assume $\Phi^{-1}(G \cdot
\alpha)$ for $\alpha \in \R_+$ consists of a single $G$-orbit. Then
every neighborhood of $\alpha$ in $\R_+$ contains a smaller
neighborhood $V$ such that the quotient $\Phi^{-1}(G \cdot V)/G$ is
contractible. 

Moreover, any reduced space $\Phi^{-1}(G \cdot \beta)/G$ of
complexity one is homeomorphic to a $2$-sphere 
for $\beta \in V$ and $\beta \neq \alpha$.
\end{lemma}

\begin{proof}
Assume $\alpha =0$, and $\Phi^{-1}(0)$ consists of a single orbit.
The stabilizer of this single orbit is either $S^1$, or a finite
subgroup $\Gamma \subset G$. The local model of the orbit  
is $G\times_{S^1} (\R^2 \times \C)$ or $G \times_\Gamma \R^3$.  
The quotients $(G \times_{S^1} (\R^2\times \C))/G$ and 
$(G \times_\Gamma \R^3)/G$ 
are homeomorphic to $(\R^2 \times\C)/S^1$  and 
$\R^3/\Gamma$ respectively. 
They are both contractible. For $\beta \neq 0$ in the local
model, direct computation shows that 
$\Phi^{-1}(G \cdot \beta)/G = \Phi^{-1}(\beta)/S^1$ is  
homeomorphic to $S^2$. 

For $\alpha \neq 0$, it suffices to show that the statements hold on
the symplectic cross-section $X$, which is a complexity one 
Hamiltonian $S^1$-manifold. It is immediate from \cite{kt1}.
\end{proof}

\comment{need to explain that local model is the neighborhood of the
short fiber?}

\begin{lemma}\labell{L:sbundleM}
Let $G$ be $\Su(2)$ or $\So(3)$ and 
let $(M, \ow, \Phi)$ be a six dimensional Hamiltonian $G$-manifold
with a moment map $\Phi\colon M \to \R^3$. Denote by $M^o$ the union of all
the tall moment fibers in $M$ and let $I^o = \Phi(M^o) \cap \R_+$. Then for
every $\alpha \in I^o$, there exists a neighborhood $V \subset \R_+$
such that $\Phi^{-1}(G \cdot V)$ is homeomorphic to 
$\Sigma \times (V \cap I^o)$, where $\Sigma$ is a surface. 
In other words, $\overline{\Phi}\colon  M^o/G \to I^o$ is topologically a
surface bundle where the map $\overline{\Phi}$ is induced by the
moment map. 
\end{lemma}

\comment{this is corollary to flattening}

\comment{need to explain $\overline{\Phi}$ somewhere}

\begin{proof}
If $\alpha =0$, this is Corollary \ref{C:sbundle}. If
$\alpha \neq 0$, we deduce from Proposition \ref{P:alphaflat}.
\end{proof}

In particular, the set $I^o$ of points in $\R_+$ whose moment fiber
is tall is connected. We have the following result:

\begin{corollary} 
Let $G$ be $\Su(2)$ or $\So(3)$ and 
let $(M, \ow, \Phi)$ be a six dimensional Hamiltonian $G$-manifold
with a moment map $\Phi\colon M \to \R^3$. Then all the tall symplectic
quotients have the same genus. 

Define the genus of a point to be zero. 
Then the genus is a well-defined notion for the complexity one
Hamiltonian $G$-manifold $M$. 
\end{corollary}

\begin{proof}[proof of proposition \ref{P:techM}]
We use Leray-Serre spectral sequence for the
induced map
$\overline{\Phi}\colon  M/G \to I = \Phi(M) \cap \R_+$.
Cover $I$ with open sets $\{U_i\}= \mathcal{U}$, 
then $\overline{\Phi}^{-1}\mathcal{U}$ is a cover for $M/G$. 
There is a spectral sequence converging to $H^*(M/G)$ with $E_2$ term
\[
E_2^{p,q} = H^p(\mathcal{U}, \mathcal{H}^q)
\]
where $\mathcal{H}^q(U) = H^q({\overline{\Phi}}^{-1}(U))$ 
is the presheaf on $I$.

Assume first that all the moment fibers in $M$ are tall. Then the
quotient $M/G$
by Lemma \ref{L:sbundleM} is topologically a surface bundle
over $I$. Then $\mathcal{H}^q$ is locally constant on $\mathcal{U}$
and the groups $\mathcal{H}^q(U) = 
H^q({\overline{\Phi}}^{-1}(U)) = H^q(\Sigma)$ are constant on
contractible open sets $U$. So $E_2^{0, q} = H^q(\Sigma)$, and $E_2^{p,q} =
0$ for $p \neq 0$. The restriction map $H^i(M/G; \Z) \to 
H^i(\Sigma; \Z)$ is an isomorphism.

If there exists a short moment fiber, by Lemma \ref{L:genus0}, the
genus must be $0$. Since the moment fiber is tall for all interior
points of $I$, the short fiber takes place at the end point of $I$.
We cover $I = [\alpha, \beta]$ with three connected open sets such that 
$U_0$ covers $\alpha$, $U_2$ covers $\beta$, and $U_{01},
U_{12}$ are connected and $U_{02}$ is empty. Then 
$E_1^{0, 0} = \Z \oplus \Z \oplus \Z$, $E_1^{1,0} = \Z \oplus \Z$,
$E_1^{1,2} = \Z \oplus \Z$, 
$E_1^{0,2} = F_\alpha \oplus \Z \oplus F_\beta$ where $F_* = \Z$ if
$\Phi^{-1}(*)$ is tall, and $F_* = 0$ if $\Phi^{-1}(*)$ is short, and 
all other $E_1^{p, q}$ vanish. So we have $E_2^{0,0} = \Z$,
$E_2^{p,0} = 0$ for $p \geq 1$, $E_2^{p,1} = 0$ for all $p$, and
$E_2^{0,2} = 0$. We see $E_2^{0,2} = E_2^{1,1} = E_2^{2,0} = 0$, and
$E_2^{0,1} = E_2^{1,0}=0$. So 
$H^i(M/G; \Z) = 0$ for $i=1, 2$. So the restriction to 
$H^i(\Sigma; \Z)$ is one-to-one for $i=1, 2$.    
\end{proof}

\comment{need to explain the cover $\mathcal{U}$}

\section{Passing to $M/G$}

\begin{proposition}\labell{P:eliminatewM}
Let $G$ be $\Su(2)$ or $\So(3)$ and let
$(M, \ow, \Phi)$ and $(M', \ow', \Phi')$ be compact connected
six dimensional Hamiltonian $G$-manifolds such that 
$\Phi(M) = \Phi'(M')$. There
exists an equivariant symplectomorphism from $M$ to $M'$ that
respects the moment maps if and only
if they have the same Duistermaat-Heckman function and there exists 
a $\Phi$-$G$-diffeomorphism from $M$ to $M'$. 
\end{proposition}

\begin{proof}
%The fact that these conditions are necessary is clear. 
Assume $f\colon M \to M'$ is a $\Phi$-$G$-diffeomorphism. By 
Proposition \ref{P:techM}, the restriction map 
$H^2(M/G; \Z) \to H^2(\Phi^{-1}(G\cdot\alpha)/G;\Z)$ 
is injective. The same arguments as in Section \ref{S:eliminatew} 
apply here; we only need to show that 
the 2-form $\ow_t = (1-t)\ow + t f^*\ow'$ is nondegenerate 
everywhere for $0 \leq t \leq 1$. 

Over zero, this is true by Lemma \ref{L:nondeg}. Away from zero, we
reconstruct the two-form $\ow_t$ by (\ref{E:coupling}). Lemma~3.6 in
\cite{kt1} guarantees the nondegeneracy on the symplectic 
cross-section and therefore on the manifold.  
\end{proof}

%\begin{lemma}[\cite{kt1} Lemma 3.6]\labell{L:alphanondeg}
%Let $S^1$ act effectively on a four dimensional manifold $M$. Let 
%$\ow_0$ and $\ow_1$ be $S^1$-invariant symplectic forms on $M$ with 
%the same moment map $\Phi$ and assume that $\ow_0$ and $\ow_1$ induce 
%the same orientation. Then the $S^1$-invariant two-form 
%$\ow_t=(1-t)\ow_0 + t\ow_1$ is nondegenerate for all $0\leq t \leq 1$.
%\end{lemma}

\begin{proposition}\labell{P:pass}
Let $G$ be $\Su(2)$ or $\So(3)$. 
Let $(M, \ow, \Phi)$ and $(M', \ow', \Phi')$ be compact connected
six dimensional Hamiltonian $G$-manifolds 
such that $0 \in \Phi(M) = \Phi'(M')$. 
Assume that $M$ and $M'$ have the same Duistermaat-Heckman function, 
that their principal stabilizers of the zero fibers are the same,
and that their first Stiefel-Whitney classes of the zero fibers
are the same when applicable. There
exists a $\Phi$-$G$-diffeomorphism from $M$ to $M'$ if and only
if there exists a $\Phi$-diffeomorphism from $M/G$ to $M'/G$. 
\end{proposition}

\begin{proof}
Assume there exists a $\Phi$-diffeomorphism $\Psi$ from $M/G$ to $M'/G$. 
Let $X$ denote the symplectic cross-section $\Phi^{-1}(\R_{>0})$ and
$\Phi_X$ denote its moment map.  
By Lemma \ref{L:phi121}, there exists a $\Phi_X$-diffeomorphism
$\psi$ on $X$. By \cite{kt1}, $\psi$
locally lifts to $\Phi_X$-$S^1$-diffeomorphisms. 

Pick an open $S^1$-invariant cover $\{U_i\}$ on $\Phi^{-1}(\R_+)$ such 
that $\Psi$ lifts to a $\Phi$-$G$-diffeomorphism $\Psi_1$ on $U_1$ and such
that $U_i \cap \Phi^{-1}(0) = \emptyset$ for all $i \neq 1$. 
We can further refine $\mathcal{U}$ so that 
$(U_i \cap U_j)/S^1$ is simply connected for all $i \neq j$ 
and on each $U_i$, $i\neq 1$, there exists a
$\Phi_X$-$S^1$-diffeomorphism $\psi_i\colon U_i \to 
M'$ that is a lift of $\psi$. By Lemma \ref{L:lift121}, $\Psi_1$
induces a $\Phi_X$-$S^1$-diffeomorphism $\psi_1$ 
on $U_1 \smallsetminus \Phi^{-1}(0)$ and $\psi_1$ is a lift of $\psi$.  

By Theorem \ref{T:hs} below, there
exist smooth $S^1$-invariant functions $g_{ij}\colon U_i\cap U_j \to
S^1$ such that 
$\psi_j^{-1} \circ \psi_i(m) = g_{ij}(m) \cdot m$
for all $m \in U_i \cap U_j$. These
functions form a \v{C}ech cocycle $g \in \check{C}^1(\mathcal{U},
S^1)$. If there exists $g_i\colon  U_i \to S^1$ such that $g_{ij} = 
g_i^{-1}\cdot g_j$, then $g_i \cdot \psi_i = g_j \cdot \psi_j$ on
$U_i \cap U_j$.  
Namely, $\{g_i \cdot \psi_i\}$ form a global
$\Phi_X$-$S^1$-diffeomorphism.

Consider the sheaves of $S^1$-invariant functions in $\R$, $S^1$, and
$\Z$. 
Since $(U_i \cap U_j)/{S^1}$ are simply connected, 
$H^1(U_i\cap U_j; \Z)=0$, and the exponential map
\[
\exp\colon  \R(U_i \cap U_j) \to S^1(U_i \cap U_j)
\] is surjective. 
The short exact sequence $0 \to \Z \to \R \to S^1 \to 1$ induces an
exact sequence $\check{H}^1(\mathcal{U}; \R) \to
\check{H}^1(\mathcal{U}; S^1) \to \check{H}^2(\mathcal{U}; \Z) \to 
\check{H}^2(\mathcal{U}; \R)$. Since there exists a partition of
unity on $\mathcal{U}$, the cohomology 
$\check{H}^1(\mathcal{U}; \R) = \check{H}^2(\mathcal{U}; \R)= 0$ and
hence the Bockstein operator 
$\delta\colon  \check{H}^1(\mathcal{U}; S^1)\to
\check{H}^2(\mathcal{U}; \Z)$ is an isomorphism. 

Now consider the image $[c]$ of $[g]$ in 
$\check{H}^2(\mathcal{U}; \Z)$.
Since $(U_i \cap U_j)/S^1$ is simply connected, we choose any particular
branch of the logarithm, and obtain that
$c_{ijk} =
\log{g_{jk}}-\log{g_{ik}}+\log{g_{ij}}$.  
Since $M$ and $M'$ have the same Duistermaat-Heckman function, 
\cite{kt1} asserts that there exist $b_{ij}$ such that
$c_{ijk} = b_{jk} - b_{ik} + b_{ij}$. 
We can choose a different branch of the
logarithm $\log^{\text{new}}g_{ij} =
\log^{\text{old}}g_{ij} - b_{ij}$ so that 
$c_{ijk} =
\log{g_{jk}}-\log{g_{ik}}+\log{g_{ij}}= 0$.
Take a partition of unity $\lambda_i$ subordinate to $\{U_i\}$, and
define \[
g_i = \exp\left(-\sum_k \lambda_k \log{g_{ik}}\right).
\] Then $g_i^{-1}\cdot g_j = 
\exp \left(\sum_k \lambda_k \log{g_{ik}}
-\sum_k \lambda_k \log{g_{jk}} \right) 
= \exp \left( \sum_k \lambda_k \log{g_{ij}}\right) = g_{ij}$ on
$U_i \cap U_j$. 
Note $g_1 = \exp\left(-\sum_k \lambda_k \log{g_{1k}}\right) = 1$ in $U_1
\smallsetminus U_k$ for $k \neq 1$. In particular, $g_1 = 1$ near
$\Phi^{-1}(0)$. By Lemma \ref{L:phig121}, we extend $g_i\cdot
\psi_i$ from the symplectic 
cross-section to the entire manifold to obtain a
global $\Phi$-$G$-diffeomorphism and $g_1 \cdot \Psi_1 = \Psi_1$ in a
neighborhood of $\Phi^{-1}(0)$.   
\end{proof}

\comment{choice of $U_i$? not correct proof}

%\begin{proposition}[\cite{kt1} Proposition 4.2]\labell{P:alphalift}
%Let $(M, \ow, \Phi)$ and $(M', \ow', \Phi')$ be four dimensional 
%Hamiltonian $S^1$-manifolds such that $\Phi(M)=\Phi'(M') = U$.
%Assume that for all $\alpha \in U$, the $\alpha$ fiber is tall
%and the restriction map 
%$H^2(M/S^1;\Z) \to H^2(\Phi^{-1}(\alpha)/S^1;\Z)$ is injective. 
%Then every $\Phi$-diffeomorphism from $M/S^1$ to $M'/S^1$ lifts to a 
%$\Phi$-$S^1$-diffeomorphism from $M$ to $M'$ if and only if 
%$M$ and $M'$ have the same Duistermaat-Heckman function. 
%\end{proposition}

%\comment{did i quote the wrong thing? i only need local lifts}

\begin{theorem}[\cite{hs}]\labell{T:hs}
Let $S^1$ act on a manifold $M$. Let $f\colon M\to M$ be an equivariant 
diffeomorphism that preserves the orbits. There exists a smooth
invariant function $h\colon M \to S^1$ such that $f(m) = h(m) \cdot m$ 
for all $m\in M$.
\end{theorem}

\section{Global Uniqueness}\labell{S:part2end}

Let $G$ be $\Su(2)$ or $\So(3)$ and let $(M, \ow, \Phi)$ be a 
six dimensional Hamiltonian $G$-manifold. 
Let $E$ denote the set of exceptional orbits in $M$. We consider 
the projections $M\to M/G$ and $\fg^* \to \fg^*/G$ and the 
map $\overline{\Phi}$ induced by the moment map $\Phi$.
The
\textbf{isotropy skeleton} is the space $E/G$ where each point is 
labeled by its isotropy representation, together with the map 
$\overline{\Phi}\colon  E/G \to \fg^*/G$. Two isotropy skeletons are 
considered the same if there exists a homeomorphism $f\colon  E/G \to E'/G$
that sends each point to a point with the same isotropy representation
and such that $\overline{\Phi} = \overline{\Phi'}\circ f$. 

\comment{this definition should check against local statements?}

We have the following global uniqueness theorem:

\begin{theoremb}
Let $G$ be $\Su(2)$ or $\So(3)$. 
Let $(M, \ow, \Phi)$ and $(M', \ow',\Phi')$ be compact connected
six dimensional Hamiltonian $G$-manifolds such that 
$\Phi(M)=\Phi'(M)$. Then $M$ and $M'$ are isomorphic if and only if
they have the same Duistermaat-Heckman function,
the same genus, the same isotropy skeleton,  
the same principal stabilizers of the manifolds and of 
the zero fibers, 
and the same first Stiefel-Whitney class of the zero fibers when
applicable.  
\end{theoremb}

We now introduce the final ingredient in the proof of the above theorem:

\comment{rewrite carefully}

\begin{lemma}\labell{L:patch}
Let $G$ be $\Su(2)$ or $\So(3)$ and $(M, \ow, \Phi)$ and $(M', \ow',
\Phi')$ be six dimensional Hamiltonian $G$-manifolds. 
Let $I_1$ and $I_2$ be open intervals in $\R_{>0} \subset \fg^*$
such that $I_1 \cap I_2 \neq \emptyset$. 
Let $U_i = \Phi^{-1}(G \cdot I_i)$ and $U_i'={\Phi'}^{-1}(G\cdot I_i)$ 
for $i =1, 2$.  
Assume that $U_i$ and $U_i'$ admit flattenings and that 
$U_1 \cap U_2$ contains only orbits with finite stabilizers. 
Let
$g_1\colon  U_1/G \to U'_1/G$ and $g_2\colon  U_2/G\to U_2'/G$ be
$\Phi$-diffeomorphisms. Then 
there exists a $\Phi$-diffeomorphism 
$\tilde{g}\colon (U_1 \cup U_2)/G \to (U_1'\cup U_2')/G$ such that 
$\tilde{g} = g_1$ 
on $U_1 \smallsetminus U_2$.
\end{lemma}

\begin{proof}
Without loss of generality, 
let $I_1 = (\alpha, b), \, I_2=(a, \beta)$, and 
$I_1 \cap I_2 = (a, b)$. 
First assume that $U_1 \cap U_2$ has no exceptional orbits.
Hence $(U_1 \cap U_2)/G$ is diffeomorphic to the product surface bundle
$\Sigma \times (a, b)$ where $\Sigma$ is the reduced space $
\Phi^{-1}(G\cdot\mu)/G$ at any $\mu \in \Phi(U_1 \cap U_2)$.
Similarly $(U_1' \cap U_2')/G \simeq
\Sigma' \times (a, b)$ where $\Sigma'$ is the reduced space at
any $\nu \in \Phi'(U_1' \cap U_2')$. 
Denote $\Sigma \times \{t\}$ by $\Sigma_t$ and $\Sigma' \times \{t\}$
by $\Sigma'_t$. We then set $(g_i)_t = g_i |_{\Sigma_t}$ so that 
$(g_i)_t$ is a diffeomorphism from $\Sigma_t$ to $\Sigma'_t$ for
$t \in (a, b)$. We can choose connections on these two surface
bundles so that the flows 
$f^s\colon \Sigma_t \to \Sigma_{t+s}$ and ${f'}^s\colon \Sigma_t' \to
\Sigma_{t+s}'$ of the horizontal lift of the vector field 
$\frac{\partial}{\partial t}$ on $(a, b)$ satisfy
${f'}^s \circ (g_1)_t = (g_1)_{t+s} \circ f^s$ for $t,\, t+s \in (a,
a+\epsilon)$ and ${f'}^s \circ (g_2)_t = (g_2)_{t+s} \circ f^s$
for $t, \, t+s \in (b-\epsilon, b)$ for some small $\epsilon >0$. 
Moreover, we may assume that these are also defined over $a$ and $b$.
We define $\gamma_t\colon \Sigma_t \to \Sigma'_t$ for $t \in (a, b)$ by
\begin{equation}
\gamma_t = {f'}^{t-a} \circ (g_1)_a \circ f^{-(t-a)}.
\end{equation} 
It's
easy to see that $\gamma_t = (g_1)_t$ for $t \in (a, a+\epsilon)$ and
\begin{equation}\labell{E:gammat}
\gamma_t = (g_2)_t \circ h_t \qquad \text{for} \quad t \in (b-\epsilon, b) 
\end{equation}
where $h_t \colon \Sigma_t \to \Sigma_t$
is a diffeomorphism determined 
by $(g_1)_a$, $(g_2)_b$, and the flows $f^s$ and
${f'}^s$. 

With the help of a smooth function $\rho\colon
(b-\epsilon, b) \to (b-\epsilon, b)$ such that $\rho(t)
= t$ near $b-\epsilon$ 
and $\rho(t)=1$ near $b$, we can reparametrize $h_t$ and $\gamma_t$ for 
$t \in (b-\epsilon, b)$ so that 
$\gamma_t$ in (\ref{E:gammat}) 
becomes $(g_2)_t \circ h$ for $t \in (c, b)$, where $h$ denotes the
diffeomorphism $h_b$ and $c$ is some value in $(b-\epsilon, b)$ near
$b$.  

The flattening of $U_2$ determines a 
new trivialization $\Phi^{-1}(G \cdot (c, b))/G \simeq \Sigma
\times (c, b)$. The exceptional sheets defined by this
trivialization correspond to a set of points $\{x_i\}$ on $\Sigma$. 
There exists a rigid diffeomorphism
$\lambda\colon \Sigma \to \Sigma$ such that $\lambda(h(x_i)) = x_i$. 
In fact, there exists an isotopy $\lambda_t \colon \Sigma_t \to
\Sigma_t$ such that $\lambda_t = \text{id}$ near $c$ and $\lambda_t =
\lambda$ near $b$. So we can construct a 
new map $g\colon (U_1 \cap U_2)/G \to (U_1' \cap
U_2')/G$ such that 
\begin{equation}\labell{E:newg}
g\bigr|_{\Sigma_t} = \begin{cases}
	(g_1)_t, \qquad & a < t < a+\epsilon \\
	\gamma_t, & a+\epsilon \leq t \leq b-\epsilon\\
	(g_2)_t \circ H_t,&b-\epsilon \leq t \leq b-\delta\\
	(g_2)_t \circ H, \qquad & b - \delta < t < b \\
	\end{cases}
\end{equation}
where $\epsilon > \delta >0$ are small and 
$H\colon  \Sigma \to \Sigma$ is a rigid diffeomorphism which sends
the points on $\Sigma$ that correspond to the exceptional sheets on
$U_2$ to themselves. 

By Proposition \ref{P:alpharigidext}, the rigid map
$H$ extends to a $\Phi$-diffeomorphism 
$\Tilde{H}$ from $({U_2 \smallsetminus U_1})/{G}$ to 
$({U_2 \smallsetminus U_1})/{G}$.
And the
new $\Phi$-diffeomorphism $\tilde{g}\colon (U_1 \cup U_2)/G \to (U_1'
\cup U_2')/G$ can be defined as $g_1$ on 
$(U_1 \smallsetminus 
U_2)/G$, $g$ on $(U_1 \cap U_2)/G$, and $g_2 \circ \Tilde{H}$ on $(U_2
\smallsetminus U_1)/G$. 

Now assume there exist exceptional orbits in $U_1 \cap U_2$. By
assumption, every orbit in $U_1 \cap U_2$ has
a finite stabilizer. The moment map $\Phi\colon U_1 \cap U_2 \to
\R^3$ is then a proper submersion. Since $0 \notin \Phi(U_1 \cap
U_2)$, the norm square of the moment map $|\Phi|^2$ from $U_1 \cap U_2$
to $I_1 \cap I_2$ is also a proper submersion. Therefore $U_1 \cap U_2$ is
$G$-equivariantly diffeomorphic to a trivial bundle 
$Z \times (a, b)$ 
where $Z/G$ is homeomorphic to the reduced space 
$\Phi^{-1}(G\cdot \mu)/G \simeq \Sigma$. 

Let $g_i\colon U_i/G \to U_i'/G$ lift to $\Phi$-$G$-diffeomorphisms
$\Hat{g}_i\colon U_i \to U_i'$. 
We proceed as before, but take $G$-equivariant connections
consistent with $\Hat{g}_i$. Again we use the smooth function
$\rho$ to reparametrize so that we obtain a new map $\Hat{g}\colon 
U_1 \cap U_2 \to U_1'\cap U_2'$ such that
\[
\Hat{g}\bigr|_{Z_t} = \begin{cases}
	(\Hat{g}_1)_t, \qquad & a < t < a+\epsilon \\
	\Hat{\gamma}_t, & a+\epsilon \leq t \leq b-\epsilon\\
	(\Hat{g}_2)_t \circ \Hat{h}_t,&b-\epsilon \leq t \leq b-\delta\\
	(\Hat{g}_2)_t \circ \Hat{h}, \qquad & b - \delta < t < b \\
	\end{cases}
\]   
We then return to the $\Phi$-diffeomorphism level. Since 
a $\Phi$-diffeomorphism induced from a $\Phi$-$G$-diffeomorphism
takes each exceptional orbit to an exceptional orbit with
the same stabilizer, we only need to construct 
a rigid diffeomorphism that sends the
exceptional sheets defined by the flattening of $U_2$ to themselves. 
This and the rest follow the same arguments as before. 
\end{proof}

\comment{exceptional orbit assumption?}

%\begin{proposition}[\cite{kt1} Proposition 11.1]\labell{P:alpharigidext}
%Let $(M, \ow, \Phi)$ and $(M', \ow', \Phi')$ be four dimensional
%Hamiltonian $S^1$-manifolds such that $\Phi(M)=\Phi'(M')=U$. Assume that 
%$M$ and $M'$ contain only tall fibers and that they admit flattenings about 
%$\alpha \in U$. Let $\Sigma$ and $\Sigma'$ denote the associated marked
%surfaces. Then any rigid map $h\colon \Sigma \to \Sigma'$ extends to a 
%$\Phi$-diffeomorphism $g\colon M/S^1 \to M'/S^1$.
%\end{proposition}

We also need the following lemma: (See \cite{k1} Appendix B, and \cite{s}.)

\begin{lemma} \labell{L:smale}
Let $k\colon S^2 \to S^2$ be a diffeomorphism such that $k$ is
a rotation on a neighborhood of the north pole and on a neighborhood
of the south pole. Then there exists an isotopy $k_t\colon 
S^2 \to S^2$,  with $k_0 = k$ and $k_1 = $ identity, such that each
$k_t$ is a rotation on a neighborhood of each pole. 
\end{lemma}

\mute{smale}{
\begin{lemma}\labell{L:patch2}
Let $G$ be $\Su(2)$ or $\So(3)$ and $(M, \ow, \Phi)$ and $(M', \ow',
\Phi')$ be six dimensional Hamiltonian $G$-manifolds. 
Let $I_1$ and $I_2$ be open intervals in $\R_{>0} \subset \fg^*$
such that $I_1 \cap I_2 \neq \emptyset$. 
Let $U_i = \Phi^{-1}(G \cdot I_i)$ and $U_i'={\Phi'}^{-1}(G\cdot I_i)$ 
for $i =1, 2$.  
Assume that $U_i$ and $U_i'$ admit flattenings and that 
$U_1 \cap U_2$ contains only orbits with finite stabilizers. 
Let $g_1\colon  U_1/G \to U'_1/G$ and $g_2\colon  U_2/G\to U_2'/G$ be
$\Phi$-diffeomorphisms. Then 
there exists a $\Phi$-diffeomorphism 
$\tilde{g}\colon (U_1 \cup U_2)/G \to (U_1'\cup U_2')/G$ such that 
$\tilde{g} = g_1$ 
on $U_1 \smallsetminus U_2$.
\end{lemma}

\comment{Need to distinguish $U, U/G, \Phi(U)$, etc and whether we
are on the symplectic cross-sections or not}}

\begin{proof}[Proof of Theorem B]
%The steps are:
%\begin{enumerate}
%\item By the invariants, there exist local $\Phi$-diffeomorphisms.
%\item WLOG, can have a nice covering. 
%There exists a cover $\{I_j\}$ of $\Phi(M) \cap \R_+$ such that
%$g_j\colon \Phi^{-1}(G\cdot I_j)/G \to {\Phi'}^{-1}(G \cdot I_j)/G$ are
%$\Phi$-diffeomorphisms. We can assume that 
%$\min I_j < \min I_{j+1}$ and that $U_j$ have flattenings if they
%contain only tall fibers. 
%\item Exclude $0$. Use induction on $g_1, \dots, g_{n-1}$.
%\item If $\Sigma \neq S^2$, then use induction on $g_n$.
%\item If $\Sigma \simeq S^2$, then use Smale's Theorem and make
%$g=g_n$. 
%\item So we have a global $\Phi$-diffeomorphism. So we have a global
%$\Phi$-$G$-diffeomorphism. So we have a global isomorphism.  
%\end{enumerate} 
If $0 \notin \Phi(M)$, this is Proposition \ref{P:nonzero}. Assume 
the moment image contains zero. 
By the assumption on the genus and the isotropy skeletons, 
there exist a $G$-invariant neighborhood $V$ of 
$\alpha$ such that $\Phi^{-1}(V)/G$ and $\Phi^{-1}(V)/G$ are 
$\Phi$-diffeomorphic for each $\alpha \in \fg^*$. 
We can cover 
$\Phi(M) \cap \R_+$ by connected open intervals such that
$g_j\colon \Phi^{-1}(G\cdot I_j)/G \to {\Phi'}^{-1}(G \cdot I_j)/G$ are
$\Phi$-diffeomorphisms. Taking refinements if necessary, we can assume that 
$\{I_j\}$ have no triple intersections, that $\min I_j < \min I_{j+1}$,
that $U_j = \Phi^{-1}(G \cdot I_j)$ have flattenings if they
contain only tall fibers, and that $U_j \cap U_{j+1}$ contain only orbits with
finite stabilizers. 

Replace $U_1$ by $U_1 \smallsetminus \Phi^{-1}(0)$. We can use induction 
and Lemma \ref{L:patch} to modify $g_1, \dots, g_{n-1}$ so that we obtain    
a new $\Phi$-diffeomorphism $\Tilde{g}_1\colon (U_1\cup \cdots \cup U_{n-1})/G
\to (U_1' \cup \cdots \cup U_{n-1}')/G$ with $\Tilde{g}_1 = g_1$ on 
$U_1 \smallsetminus U_2$. 

If the reduced space $\Phi^{-1}(G\cdot\alpha)/G$ is tall for all 
$\alpha \in I_n$, we can use induction 
again on $U_n$ so that we obtain a new $\Phi$-diffeomorphism 
$\Tilde{g}\colon (U_1\cup \cdots \cup U_{n})/G
\to (U_1' \cup \cdots \cup U_{n}')/G$ with $\Tilde{g} = g_1$ on 
$U_1 \smallsetminus U_2$. So $\Tilde{g}$ naturally extends to 
$\Phi^{-1}(0)/G$ and we obtain a global $\Phi$-diffeomorphism 
$\Tilde{g}\colon M/G \to M'/G$. 

By Lemmas \ref{L:genus0} and \ref{L:sbundleM},  
if the reduced space $\Phi^{-1}(G\cdot\beta)/G$ is short for $\beta = 
\max I_n$, the reduced space $\Phi^{-1}(G\cdot\alpha)/G$ is homeomorphic 
to $S^2$ for any $\alpha \neq \beta$ in $I_n$. 
By \cite{k1}, there exist at most two exceptional orbits in 
$\Phi^{-1}(G\cdot\alpha)$. 

Replace $I_n$ by $I_n \smallsetminus \{\beta\}$ and $U_n$ by $U_n 
\smallsetminus \Phi^{-1}(G\cdot\beta)$. 
Denote $I_{n-1}\cap I_n$ by $(a, b)$. We construct a map 
$g\colon (U_{n-1}\cap U_n)/G \to (U_{n-1}' \cap U_n')/G$  
as in the proof of Lemma \ref{L:patch} up to the form in (\ref{E:newg}).
By Lemma \ref{L:smale}, there exists a rigid isotopy 
$k_t\colon S^2 \to S^2$ such that $k_0 = H$ and 
$k_1 = \text{identity}$. So we obtain a new map 
$g\colon (U_{n-1} \cap U_n)/G \to (U_{n-1}' \cap U_n')/G$ such that
\[
g\bigr|_{\Sigma_t} = \begin{cases}
	(\Tilde{g}_1)_t, \qquad & a < t < a+\epsilon \\
	\gamma_t, & a+\epsilon \leq t \leq b-\epsilon\\
	(g_n)_t \circ H_t,&b-\epsilon \leq t \leq b-\delta\\
%	(g_n)_t \circ H, \qquad & b - \delta \leq t \leq b - \delta/2 \\
	(g_n)_t, \qquad & b - \delta < t < b\\
	\end{cases}
\]
where $\epsilon > \delta >0$ are small numbers. 
%and $H\colon  \Sigma \to \Sigma$ is a rigid diffeomorphism which sends
%the points on $\Sigma$ that correspond to the exceptional sheets on
%$U_n$ to themselves. 
Therefore the new $\Phi$-diffeomorphism $\Tilde{g}$ given by
$\Tilde{g}_1$ on $((U_1 \cup \cdots \cup U_{n-1}) \smallsetminus U_n)/G$,
$g$ on $((U_1 \cup \cdots \cup U_{n-1}) \cap U_n)/G$, and $g_n$ on 
$(U_n \smallsetminus (U_1 \cup \cdots \cup U_{n-1}))/G$ is well-defined 
on $U_1 \cup \cdots \cup U_n$ such that $\Tilde{g} = g_1$ on 
$U_1 \smallsetminus U_2$ and $\Tilde{g} = g_n$ on $U_n \smallsetminus U_{n-1}$.
Hence $\Tilde{g}$ naturally extends to $\Phi^{-1}(0)/G$ as well as 
$\Phi^{-1}(G\cdot \alpha)/G$, where $\alpha = \max I_n$.
The map $\Tilde{g} \colon M/G \to M'/G$ is a global 
$\Phi$-diffeomorphism. 

By Proposition \ref{P:pass}, there exists a $\Phi$-$G$-diffeomorphism from 
$M$ to $M'$. By Proposition \ref{P:eliminatewM}, there exists an equivariant
symplectomorphism from $M$ to $M'$.
\end{proof}

\comment{need to check the reference}


\begin{thebibliography}{ABCDE}
	\bibitem[AH]{ah}
		K.~Ahara and A.~Hattori, 
		\emph{4 dimensional symplectic $S^1$-manifolds
		admitting moment map},
		J. Fac. Sci. Univ. Tokyo Sect. IA, Math.\textbf{38} (1991),
		251-298.
	\bibitem[Au1]{aU1}
		M.~Audin,
		\emph{Hamiltoniens p\'eriodiques sur les vari\'et\'es
		symplectiques 
		compactes de dimension 4},
		g\'eometrie symplectique et m\'ecanique, Proceedings
		1998, C.~Albert ed.~,
		Springer Lecture Notes in Math.\textbf{1416} (1990)
	\bibitem[Au2]{aU2}
		M.~Audin,
		\emph{The topology of torus actions on symplectic
		manifolds},
		Progress in Mathematics \textbf{93}, Birkhuser Verlag,
		basel, 1991.
%	\bibitem[BB]{bb}
%		A.~Bialynicki-Birula,
%		\emph{Remarks on the action of  an algebraic torus on
%		$k^n$}
%		II,
%		Bull. Acad. Polon. Sci. S\'er. Sci. Math. Astronom. Phys.
%		\textbf{15},
%		(1967), 123-125.
	\bibitem[BL]{bl}
		L.~Bates and E.~Lerman,
		\emph{Proper group actions and symplectic stratified spaces},
		Pacific Journal of Mathematics Vol.181, No.2, (1997), 201-229.
	\bibitem[D1]{dE1}
		T.~Delzant,
		\emph{Hamiltoniens p\'eriodiques et image convexe de 
			l'application moment}, 
		Bull. Soc. Math. France \textbf{116} (1988), 315-339.
	\bibitem[D2]{dE2}
		T.~Delzant,
		\emph{Classification des actions hamitoniennes
		compl\`ement int\'egrables de rang deux},
		Ann. Global Anal. Geom.\textbf{8} (1990), 87-112.
%	\bibitem[FK]{fk}
%		K.~-H.~Fieseler and L.~ Kaup, 
%		\emph{On the geometry of affine algebraic $\C^*$-surfaces}, 
%		Problems in the theory of surfaces and their
%		classification 
%		(Cortona, 1988), 111--140, 
%		Sympos. Math., XXXII, 
%		Academic Press, London, 1991. 
	\bibitem[F]{fI}
		R.~Fintushel, 
		\emph{Circle actions on simply connected $4$-manifolds},
		Trans. Amer. Math. Soc. \textbf{230} (1977), 147--17
	%\bibitem[GSj]{gsj}
	%	V.~Guillemin and R.~Sjamaar,
	%	lecture at Mass. Inst. Tech. , 1995.
	\bibitem[GLS]{gls}
		V.~Guillemin, E.~Lerman and S.~Sternberg,
		\emph{Symplectic fibrations and multiplicity diagrams},
		Cambridge University Press, (1995).
	\bibitem[GS1]{gs1}
		V.~Guillemin and S.~Sternberg,
		\emph{A normal form for the moment map},
		Differential geometric methods in mathmatical physics, 
		(S.~Sternberg, Ed.) Reidel, Dordrecht, Holland, 1984.
	\bibitem[GS2]{gs2}
		V.~Guillemin and S.~Sternberg,
		\emph{Symplectic techniques in physics}, 
		Cambridge University Press, (1990).
	\bibitem[HS]{hs}
		A.~Haefliger and E.~Salem,
		\emph{Actions of tori on orbifolds},
		Ann. Global Anal. Geom. \textbf{9} (1991), 37-59. 
	\bibitem[Ig]{iG}
		P.~Iglesias,
		\emph{Les $\So(3)$-vari\'et\'es symplectiques et leurs
		classification en dimension 4},
		Bull. Soc. Math. France \textbf{119} (1991), 371-396.
   	\bibitem[K1]{k1}
		Y.~Karshon,
		\emph{Periodic Hamiltonian flows on four dimensional
		manifolds},
		Mem. Amer. Math. Soc.\textbf{672} (1999).
	\bibitem[K2]{k2}
		Y.~Karshon,
		\emph{Appendix to Symplectic packings and algebraic
		geometry}
		by D.~MuDuff and L.~Polterovich,
		Invent. Math. \textbf{115}, 431-434 (1994). 
	\bibitem[KT1]{kt1}
		Y.~Karshon and S.~Tolman, 
		\emph{Centered complexity one Hamiltonian torus
		actions},
		Trans. Amer. Math. Soc. \textbf{353}  (2001), 4831-4861.
	\bibitem[KT2]{kt2}
		Y.~Karshon and S.~Tolman,
		\emph{Tall complexity one Hamiltonian torus actions},
		math.SG/0202167.
%	\bibitem[KKMS]{kkms}
%		G.~Kempf, F.~F.~Knudsen, D.~Mumford, B.~Saint-Donat,
%		\emph{Toroidal embeddings. I.}
%		Lecture Notes in Mathematics, Vol. \textbf{339}, 
%		Springer-Verlag, Berlin-New York, 1973. 
	\bibitem[Kn]{kN}
		F.~Knop, in preparation.
	\bibitem[Ko]{kO}
		A.~Kosinski, 
		\emph{Differential manifolds},
		Pure and Applied Mathematics, \textbf{138} 
		Academic Press, Inc., Boston, MA, 1993. 
	\bibitem[L]{lE}
		E.~Lerman, 
		\emph{Symplectic cuts},
		Math. Res. Lett. \textbf{2} (1995), no. 3, 247--258.
	\bibitem[LS]{ls}
		E.~Lerman and R.~Sjamaar, 
		\emph{Stratified symplectic spaces and reduction},
		Ann. of Math. (2) \textbf{134} (1991), no. 2, 375--42
%	\bibitem[LT]{lt}
%		E.~Lerman and S.~Tolman, 
%		\emph{Hamiltonian torus actions on symplectic
%		orbifolds 
%		and toric varieties},
%		Trans. Amer. Math. Soc. \textbf{349} (1997), no. 10, 
%		4201--4230. 
%	\bibitem[LV]{lv}
%		D.~Luna and Th.~Vust, 
%		\emph{Plongements d'espaces homog\`enes},
%		Comment. Math. Helv. \textbf{58} (1983), no. 2, 186--245. 
	\bibitem[M]{mA}
		C.~M.~Marle, 
		\emph{Mod\`ele d'action hamiltonienne d'un groupe de
		Lie 
		sur une vari\'et\'e symplectique},
		Rend. Sem. Mat. Univ. Politec. Torino \textbf{43}
		(1985), 
		no. 2, 227--251. 
	\bibitem[Mo]{mO}
	    J.~Moser,
		\emph{On the volume elements on a manifold},
		Trans. Amer. Math. Soc. \textbf{120} (1965), 286-294.
%	\bibitem[OR]{or}
%		P.~Orlik and F.~Raymond, 
%		\emph{Actions of the torus on $4$-manifolds. II.}
%		Topology \textbf{13} (1974), 89--112. 
	\bibitem[OW]{ow}
		P.~Orlik and P.~ Wagreich,
		\emph{Algebraic surfaces with $k^*$-action}, 
		Acta Math. \textbf{138} (1977), no. 1-2, 43--81. 
%	\bibitem[R]{rY}
%		J.~Rynes,
%		\emph{Nonsingular affine $k^*$-surfaces},
%		Trans. Amer. Math. Soc. \textbf{332} (1992), 889-921.
	\bibitem[Sch1]{sCh1}
		G.~W.~Schwarz, 
		\emph{Smooth functions invariant under the action 
		of a compact Lie group},
		Topology \textbf{14} (1975), 63--68.
	\bibitem[Sch2]{sCh2}
		G.~W.~Schwarz, 
		\emph{Lifting smooth homotopies of orbit spaces},
		Inst. Hautes \'Etudes Sci. Publ. Math. No. 51, 
		(1980), 37--135. 
	\bibitem[S]{s}
		S.~Smale,
		\emph{Diffeomorphisms of the 2-spheres},
		Proc. Amer. Math. Soc. \textbf{10} (1959), 621--626.
	\bibitem[T1]{tI1}
		D.~A.~Timash\"ev, 
		\emph{$G$-manifolds of complexity $1$},
		(Russian) Uspekhi Mat. Nauk \textbf{51} (1996), 
		no. 3(309),213--214; English translation in Russian 
		Math. Surveys \textbf{51} (1996), no. 3, 567--568.
	\bibitem[T2]{tI2}
		D.~A.~Timash\"ev,
		\emph{Classification of $G$-manifolds of complexity
		$1$},
		 (Russian) Izv. Ross. Akad. Nauk Ser. Mat. \textbf{61}
		(1997), no. 2,127--162; English translation in Izv. 
		Math. \textbf{61} (1997), no. 2, 363--397.
	\bibitem[W1]{w1}
		A.~Weinstein,
		\emph{Lectures on Symplectic Manifolds},
		CBMS Reg. conf. ser. in Math. \textbf{29}, 1977.
	\bibitem[W]{w}
		C.~Woodward,
		\emph{The classification of transversal multiplicity-free
		group actions},
		Ann. Global Anal. Geom. \textbf{14} (1996), 3-42.
\end{thebibliography}
\end{document}